\documentclass[dvipsnames]{amsart}
\usepackage[utf8]{inputenc}

\usepackage{pgf,tikz,pgfplots}
\pgfplotsset{compat=1.15}

\definecolor{gray}{rgb}{0.4,0.4,0.4}
\definecolor{wrwrwr}{rgb}{0.3803921568627451,0.3803921568627451,0.3803921568627451}
\definecolor{rvwvcq}{rgb}{0.08235294117647059,0.396078431372549,0.7529411764705882}
\definecolor{uuuuuu}{rgb}{0.26666666666666666,0.26666666666666666,0.26666666666666666}
\definecolor{red}{rgb}{0.8,0,0}
\definecolor{blue}{rgb}{0,0,1}
\definecolor{zzttff}{rgb}{0.6,0.2,1}
\usepackage{ifxetex,ifluatex}
\newif\ifxetexorluatex
\ifxetex
  \xetexorluatextrue
\else
  \ifluatex
    \xetexorluatextrue
  \else
    \xetexorluatexfalse
  \fi
\fi

\ifxetexorluatex
  \usepackage{fontspec}
\else
  \usepackage[T1]{fontenc}
  \usepackage[utf8]{inputenc}
  \usepackage{lmodern}
\fi

\usepackage[hidelinks]{hyperref}
\usepackage{amsthm}
\usepackage{amsmath}
\usepackage{amsfonts}
\usepackage{amssymb}
\usepackage{array}
\usepackage{multirow}
\usepackage{mathtools}
\usepackage[labelformat=simple, labelfont=normalfont]{subcaption}
\usepackage{comment}
\usepackage{mathrsfs}
\usepackage[margin=1in]{geometry}
\usepackage{listings}
\usepackage[capitalise, noabbrev]{cleveref}
\usepackage[shortlabels]{enumitem}
\usepackage{textcomp}

\usepackage{graphicx}
\usepackage{xcolor}
\usepackage{caption}
\usepackage{float}

\usetikzlibrary{arrows.meta, fit, calc, decorations.markings, positioning, shapes.geometric}
\tikzset{baseline={($ (current bounding box.west) - (0,1ex) $)}, auto}
\tikzset{vertex/.style={circle, inner sep=1.5pt, fill}, edge/.style={thick, line join=bevel}, optional/.style={black!50, dashed}, highlight/.style={red, very thick}, starpoint/.style={fill=white, draw=black}}

\DeclarePairedDelimiter\abs{\lvert}{\rvert}

\newcommand\set[1]{\{ #1 \}}

\renewcommand{\emptyset}{\varnothing}
\DeclareMathOperator{\link}{link}

\newtheorem{thm}{Theorem}[section]
\newtheorem{lem}[thm]{Lemma}
\newtheorem{cor}[thm]{Corollary}

\newtheorem{prop}[thm]{Proposition}
\theoremstyle{definition}
\newtheorem{defn}[thm]{Definition}

\newtheorem{rem}[thm]{Remark}
\Crefname{thm}{Theorem}{Theorems}
\Crefname{lem}{Lemma}{Lemmas}
\Crefname{cor}{Corollary}{Corollaries}
\Crefname{conj}{Conjecture}{Conjectures}
\Crefname{prop}{Proposition}{Propositions}
\Crefname{defn}{Definition}{Definitions}
\Crefname{claim}{Claim}{Claims}
\Crefname{rem}{Remark}{Remarks}
\crefname{page}{page}{pages}

\title{A characterization of two-dimensional\\ Buchsbaum matching complexes}

\author{Bennet Goeckner} 
\address{Department of Mathematics, University of San Diego}
\email{bgoeckner@sandiego.edu}

\author{Fran Herr}
\address{Department of Mathematics, University of Washington} 
\email{fanhrr@gmail.com}

\author{Legrand Jones II} 
\address{Department of Mathematics, Indiana University}
\email{legjones@iu.edu}

\author{Rowan Rowlands} 
\address{Department of Mathematics, University of Washington}
\email{rowanr@uw.edu}

\begin{document}
\maketitle

\begin{abstract}
    The matching complex $M(G)$ of a graph $G$ is the set of all matchings in $G$. A Buchsbaum simplicial complex is a generalization of both a homology manifold and a Cohen--Macaulay complex. We give a complete characterization of the graphs $G$ for which $M(G)$ is a two-dimensional Buchsbaum complex. As an intermediate step, we determine which graphs have matching complexes that are themselves connected graphs.
\end{abstract}

\section{Introduction}

Given a graph $G$, a matching is a collection of edges such that no two share a common endpoint. The matching complex $M(G)$, which is the set of all matchings in $G$, forms a simplicial complex. Matching complexes and their topology have been studied extensively; see, e.g., \cite{Jon, Wachs} for surveys of the field.

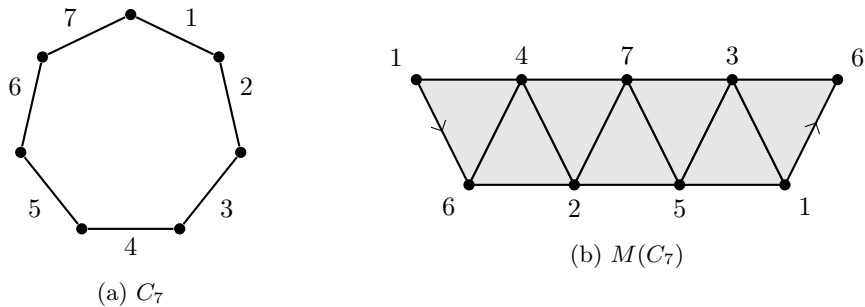
\begin{figure}
\centering
\begin{subfigure}{0.3\textwidth}
\centering
\begin{tikzpicture}

\pgftransformrotate{90}

\foreach \i in {1,...,7}
    \node (v\i) [vertex] at (360/7-\i*360/7:1.5) {};

\draw [edge] (v1) -- (v2) node [midway, auto] {1};
\draw [edge] (v2) -- (v3) node [midway, auto] {2};
\draw [edge] (v3) -- (v4) node [midway, auto] {3};
\draw [edge] (v4) -- (v5) node [midway, auto] {4};
\draw [edge] (v5) -- (v6) node [midway, auto] {5};
\draw [edge] (v6) -- (v7) node [midway, auto] {6};
\draw [edge] (v7) -- (v1) node [midway, auto] {7};

\end{tikzpicture}
\caption{$C_7$} \label{fig:Manifold-Example-C7}
\end{subfigure}
\qquad
\begin{subfigure}{0.4\textwidth}
\centering
\begin{tikzpicture}[scale=.7]
\draw[fill=black!10, edge] (1,0)--(2,2)--(0,2)--cycle;
\draw[fill=black!10, edge] (1,0)--(3,0)--(2,2)--cycle;
\draw[fill=black!10, edge] (3,0)--(4,2)--(2,2)--cycle;
\draw[fill=black!10, edge] (3,0)--(5,0)--(4,2)--cycle;
\draw[fill=black!10, edge] (5,0)--(6,2)--(4,2)--cycle;
\draw[fill=black!10, edge] (5,0)--(7,0)--(6,2)--cycle;
\draw[fill=black!10, edge] (7,0)--(8,2)--(6,2)--cycle;

\foreach \pos/\name/\labelpos in {(1,0)/6/below left, (3,0)/2/below, (5,0)/5/below, (7,0)/1/below right, (0,2)/1/above left, (2,2)/4/above, (4,2)/7/above, (6,2)/3/above, (8,2)/6/above right}
    \node [vertex, label=\labelpos:$\name$] at \pos {};

\draw [edge, decoration={markings, mark=at position 0.5 with {\arrow{Straight Barb[]}}}, decorate] (0,2)--(1,0);
\draw [edge, decoration={markings, mark=at position 0.6 with {\arrow{Straight Barb[]}}}, decorate] (7,0)--(8,2);
\end{tikzpicture}
\caption{$M(C_7)$} \label{fig:Manifold-Example-matching}
\end{subfigure}
    \caption{The matching complex of $C_7$ is a triangulated M\"{o}bius strip. Faces of $M(C_7)$ with the same label are identified.}
    \label{fig:Manifold-Example}
\end{figure}

Recently, all homology manifolds that arise as matching complexes have been classified \cite{BGM}. Outside of dimension two, all such complexes are combinatorial (i.e., PL) balls and spheres. In dimension two, more examples appear, including a torus and a M\"obius strip. See \cref{fig:Manifold-Example} for one such example.

In this paper, we characterize all graphs $G$ for which $M(G)$ is a two-dimensional Buchsbaum complex, which partially answers a question from \cite[Section~6]{BGM}. Buchsbaum complexes are a generalization of both homology manifolds and Cohen--Macaulay complexes. Though originally defined algebraically, the Buchsbaum condition is in fact a topological property \cite{Schenzel}. In dimension two, Buchsbaum complexes can be classified in terms of certain subcomplexes being connected graphs, and this is the notion we will use.

In \cref{sec:Prelim}, we introduce relevant terminology and background. In \cref{sec:Link} we classify all graphs $G$ for which $M(G)$ is a connected graph, which allows us to characterize all one-dimensional Buchsbaum and Cohen--Macaulay matching complexes in \cref{prop:CM-Buchs-1-dim}. Then we consider the local behavior of graphs with two-dimensional Buchsbaum matching complexes. \Cref{sec:Families} gives an explicit description of all graphs $G$ such that $M(G)$ is a two-dimensional Buchsbaum complex in \cref{thm:MainTheorem} and then shows that this list is exhaustive. We end with a brief discussion of similar questions in higher dimensions in \cref{sec:ConcludingRemarks}.

\section{Preliminaries}\label{sec:Prelim}

Our two main objects of study are simple graphs and simplicial complexes. For all terms not defined here, see standard references such as \cite{West} and \cite{GreenStanley}.

A \emph{(simple) graph} $G=(V,E)$ consists of a vertex set $V=V(G)$ and an edge set $E=E(G)$ whose members are two-element subsets of $V$. If $e=\set{a,b} \in E,$ we refer to vertices $a$ and $b$ as the \emph{endpoints} of the edge $e$; we will often use the notation $e=ab$. Given a graph $G$, a \emph{matching} is a collection of edges of $G$ such that no two share an endpoint. Unless stated otherwise, we will assume that all graphs mentioned in theorem statements are simple and do not have isolated vertices (i.e., vertices that are not the endpoints of any edges).

A graph is \emph{bipartite} if the vertices can be partitioned into two sets $V_1$ and $V_2$ so that every edge has one endpoint in $V_1$ and one endpoint in $V_2$, or, equivalently, if the graph contains no odd-sized cycles. We often refer to several common graphs: $K_n$ the complete graph on $n$ vertices, $C_n$ the cycle on $n$ vertices, and $S_n$ the star graph with $n+1$ vertices. We often refer to a path on $n+1$ vertices as a path of length $n$ and denote it as $P_{n+1}$.

A \emph{simplicial complex} $\Delta$ is a collection of sets with the property that if $\sigma \in \Delta$ and $\tau \subseteq \sigma$, then $\tau \in \Delta$. An element $\sigma \in \Delta$ is called a \emph{face}; throughout, we will use the convention of writing $abc$ in place of $\set{a,b,c}$, etc.\ for faces of simplicial complexes. The \emph{dimension} of a face $\sigma$ is $\dim \sigma := \abs{\sigma} - 1$, and the \emph{dimension} of $\Delta$, denoted $\dim \Delta$, is the maximum of the dimensions of its faces. A complex is \emph{pure} if all maximal faces have the same dimension. Faces of dimension $0$ and $1$ are called vertices and edges respectively.

Note that simple graphs can be thought of as simplicial complexes of dimension one (or less, if the graph has no edges). Throughout, we will often blur the distinctions between graphs and $1$-dimensional simplicial complexes and between simplicial complexes and their geometric realizations.

Given a face $\sigma \in \Delta$, its \emph{link}, denoted $\link_\Delta \sigma$ (or simply $\link \sigma$ if $\Delta$ is unambiguous), is
$$
\link_\Delta \sigma = \set{ \tau \in \Delta \mid \tau \cup \sigma \in \Delta ~\text{and}~ \tau \cap \sigma = \emptyset }.
$$
For example, the link of vertex $7$ in \cref{fig:Manifold-Example-matching} is the path with edges $42,25,53$ and the link of edge $16$ is the pair of isolated vertices $3$ and $4$. The link of a face captures the local structure of $\Delta$ near that face, and many properties of simplicial complexes---including Buchsbaumness---can be defined in terms of links.

The \emph{matching complex} $M(G)$ is the set of all matchings in $G$. Since any subset of a matching is also a matching, $M(G)$ is a simplicial complex. Note that the vertices of $M(G)$ correspond to the edges of $G$. For especially small graphs, the matching complex is easy to calculate by hand. Alternatively, we define a Mathematica function ``\lstinline|MatchingComplex|'' in the appendix.

The figures in this paper occasionally include dotted edges. If a pendant edge is dotted, then our arguments relating to that figure allow for arbitrarily many copies of that pendant. Similarly, if a figure includes a dotted path of length $2$ that only touches the rest of the graph at the ends of the path, then the graph may include arbitrarily many paths of length $2$ attached at the same points.

Many results on matching complexes concern the topological properties of their geometric realizations. Most previous results consider $M(G)$ for a family of graphs (see, e.g., \cite{Jon} and \cite{Wachs} for a survey of these results), but we will instead specify the properties of $M(G)$ and then determine the structure of $G$. Motivated by a question from \cite[Section~6]{BGM}, we will be primarily interested in complexes that satisfy the following definition.

\begin{defn}\label{def:Buchs}
Let $\Delta$ be a two-dimensional simplicial complex. We say that $\Delta$ is \emph{Buchsbaum} if for each vertex $v \in \Delta$, $\link_\Delta v$ is a connected graph with at least one edge.
\end{defn}

Note that if $\Delta$ has any maximal faces of dimension $0$ or $1$, then $\Delta$ cannot satisfy \cref{def:Buchs}, so all maximal faces must have dimension $2$, i.e., $\Delta$ is pure. For example, the link of every vertex in \cref{fig:Manifold-Example-matching} is a path on four vertices, so $M(C_7)$ is Buchsbaum.

\begin{rem}\label{rem:Buchs-gen}
In this paper we focus on Buchsbaum complexes in dimension two. In general, a Buchsbaum complex can be defined as a pure complex where the $i^{th}$ reduced homology of $\link \sigma$ is trivial for all $i < \dim \Delta - \abs \sigma$ for all nonempty faces $\sigma \in \Delta$. The Buchsbaum condition was first defined in terms of algebraic properties of the complex's associated Stanley--Reisner ring. However, the combinatorial description above is equivalent for two-dimensional complexes, and, moreover, Buchsbaumness is a topological invariant (see \cite{Schenzel,Miyazaki}).
\end{rem}

\begin{rem}
Buchsbaumness is a generalization of the \emph{Cohen--Macaulay} condition, which additionally requires that $\link_{\Delta} \emptyset$ (i.e.\ $\Delta$ itself) also has vanishing $i^{th}$ homology for all $i < \dim \Delta$. There is a related and even more restrictive class known as \emph{Gorenstein} complexes. A complete characterization of Gorenstein matching complexes is implicit in \cite[Theorem~3.1]{BGM} (via \cite[Chapter~II~Theorem~5.1]{GreenStanley}) and is proved independently in \cite[Theorem~2.1]{Ni21}.
\end{rem}

In light of these remarks, we note that if $\dim \Delta = 1$, then $\Delta$ is Buchsbaum if and only if $\Delta$ is a graph with no isolated vertices, and $\Delta$ is Cohen--Macaulay if and only if it is connected. In general, a complex $\Delta$ is Buchsbaum if and only if it is pure and $\link_\Delta v$ is Cohen--Macaulay for all vertices $v$ of $\Delta$. We will consider the one-dimensional case in \cref{sec:Link}, and it will be key to developing our results for two-dimensional complexes.

We are not aware of any overt study of Buchsbaum matching complexes in the literature, but there are some results for similar properties for certain families of graphs. For example, \cite[Theorem~15]{Garst} shows that $M(K_{m,n})$ is Cohen--Macaulay if and only if $n \ge 2m-1$, and \cite[Theorem~2.3]{Zi94} shows that this is in fact equivalent to vertex decomposability for this family of matching complexes.

\section{One-dimensional matching complexes and link behavior}\label{sec:Link}

\begin{figure}
    \centering
    \begin{subfigure}{0.4\textwidth}
        \centering
        \begin{tikzpicture}
            \node (a) [vertex, label=above:$1$
            ] at (1,1) {};
            \node (b) [vertex, starpoint, label=above:$2$
            ] at (0,0) {};
            \node (c) [vertex, label=above:$3$
            ] at (1,0) {};
            \node (d) [vertex, starpoint, label=above:$4$
            ] at (2,0) {};
            \node (e) [vertex, label=below:$5$
            ] at (1,-1) {};
            
            \draw [edge] (a)--(b)--(c)--(d)--(e);
            \draw [edge, optional] (a)--(d) (b)--(e);
        \end{tikzpicture}
        \caption{$\mathcal{G}_1$}
    \end{subfigure} \quad
    \begin{subfigure}{0.4\textwidth}
        \centering
        \begin{tikzpicture}[scale=1]
            \node (a) [vertex, label=above left:$1$
            ] at (0,1) {};
            \node (b) [vertex, label=below left:$2$
            ] at (0,-1) {};
            \node (c) [vertex, label=left:$3$
            ] at (0.8,0) {};
            \node (d) [vertex, starpoint, label=above:$4$
            ] at (2,0) {};
            \node (e) [vertex, label=above right:$5$] at (3,0) {};
            
            \draw [edge] (a)--(b)--(c)--(d)--(e);
            \draw [edge] (a)--(c);
            \draw [edge, optional] (a)--(d)--(b);
        \end{tikzpicture}
        \caption{$\mathcal{G}_2$}
    \end{subfigure} \\
    \begin{subfigure}{0.4\textwidth}
        \centering
        \begin{tikzpicture}[scale=1]
            \node (a) [vertex, label=above:$1$
            ] at (90:1) {};
            \node (b) [vertex, label=right:$2$
            ] at (18:1) {};
            \node (c) [vertex, label=below right:$3$
            ] at (306:1) {};
            \node (d) [vertex, label=below left:$4$
            ] at (234:1) {};
            \node (e) [vertex, label=left:$5$
            ] at (162:1) {};
            
            \draw [edge] (a)--(b)--(c)--(d)--(e);
            \draw [edge] (e)--(a);
            \draw [edge, optional] (a)--(c)--(e)--(b)--(d)--(a);
        \end{tikzpicture}
        \caption{$\mathcal{G}_3$}
    \end{subfigure} \quad
    \begin{subfigure}{0.4\textwidth}
        \centering
        \begin{tikzpicture}[scale=1]
            \node (a) [vertex, label=above left:$1$
            ] at (0,1) {};
            \node (b) [vertex, label=below left:$2$
            ] at (0,0) {};
            \node (c) [vertex, label=above:$3$
            ] at (1,0.5) {};
            \node (d) [vertex, label=above right:$4$
            ] at (2,1) {};
            \node (e) [vertex, label=below right:$5$
            ] at (2,0) {};
            
            \draw [edge] (a)--(b)--(c)--(d)--(e);
            \draw [edge] (a)--(c)--(e);
        \end{tikzpicture}
        \caption{Bowtie graph $B$}
    \end{subfigure}
    \caption{The families of connected graphs whose matching complexes are themselves connected graphs. Dotted edges are optional. For families $\mathcal{G}_1$ and $\mathcal{G}_2$, we may repeatedly introduce filled vertices and connect them to any (nonzero) number of unfilled vertices.}
    \label{fig:NP5k}
\end{figure}
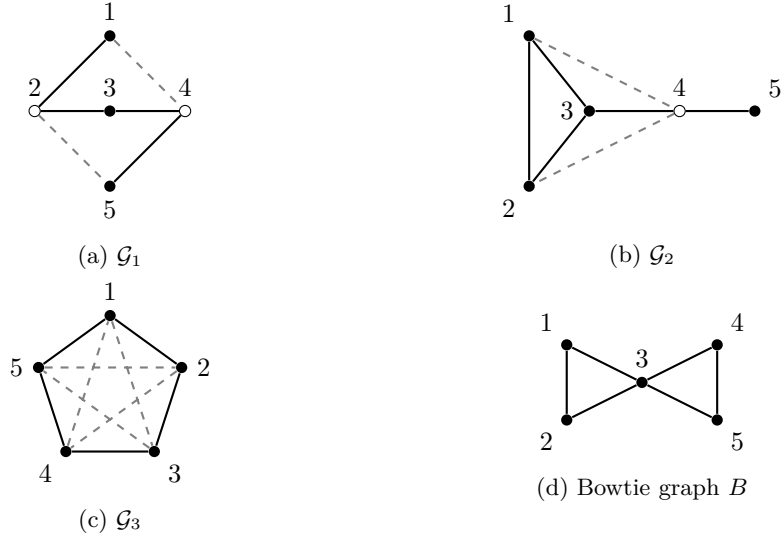

The goal of this section is to classify all Cohen--Macaulay and Buchsbaum matching complexes in dimension one. We begin by defining the families $\mathcal G_1$, $\mathcal G_2$, and $\mathcal G_3$ and the bowtie graph $B$ in \cref{fig:NP5k}. In this figure, dotted edges are optional. For families $\mathcal{G}_1$ and $\mathcal{G}_2$, we may repeatedly introduce filled vertices and connect them to any (nonzero) number of unfilled vertices.

The following is the main result in this section.

\begin{thm}\label{prop:CM-Buchs-1-dim}
Let $G$ be a graph and assume $\dim M(G)=1$.
\begin{enumerate}[(a)]
    \item $M(G)$ is Cohen--Macaulay if and only if either
    \begin{enumerate}[(i)]
    \item $G$ has two components which are each either a $K_3$ or a star graph, or
    \item $G$ is in one of the families $\mathcal{G}_1$, $\mathcal{G}_2$, and $\mathcal{G}_3$, or $G$ is the bowtie graph $B$, all of which are defined in \cref{fig:NP5k}. 
    \end{enumerate}
    
    \item $M(G)$ is Buchsbaum if and only if $G$ is a graph described in (a) or $G=K_4$ or $G=C_4$.
\end{enumerate}
\end{thm}

We will spend the rest of this section proving \cref{prop:CM-Buchs-1-dim}. We first introduce a tool that will be used throughout the remainder of the paper.

\begin{defn}
For a simple graph $G$ and an edge $e \in E(G)$, the \emph{non-adjacent subgraph} of $e$, denoted $N_e$, is the subgraph induced by all edges of $G$ that do not share any endpoints with $e$.
\end{defn}

Observe that $N_e$ will never have any isolated vertices. Furthermore, the link of the vertex $e$ in $M(G)$ is $\link_{M(G)} e = M(N_e)$. This allows us to use non-adjacent subgraphs to translate the two-dimensional Buchsbaum condition for the matching complex $M(G)$ into conditions for the graph $G$.

\begin{lem}\label{lem:2D-Buchs-Equiv}
Given a graph $G$, $M(G)$ is a two-dimensional Buchsbaum complex if and only if $M(N_e)$ is a connected graph with at least one edge for all $e \in E(G)$.
\end{lem}

\begin{proof}
A matching complex $M(G)$ is two-dimensional if and only if the largest size of a matching in $G$ is three. This is equivalent to the largest dimension of the link of a vertex in $M(G)$ being one (i.e., the link is a graph with at least one edge).

Let $M(G)$ be a two-dimensional complex. Then $M(G)$ is Buchsbaum if and only if $\link_{M(G)}e$ is a connected graph with at least one edge for each vertex $e$ of $M(G)$. Since $\link_{M(G)}e=M(N_e),$ this completes the proof.
\end{proof}

We will use the above result throughout as our main tool for characterizing two-dimensional Buchsbaum matching complexes. 

We will first consider graphs $G$ for which $M(G)$ is a connected \emph{graph}, i.e., a $1$-dimensional simplicial complex. These graphs will be instrumental in \cref{sec:Families}, and they answer the question for $1$-dimensional Cohen--Macaulay and Buchsbaum matching complexes. We first turn our attention to matching complexes of disconnected graphs.

\begin{lem}\label{lem:Disconnected-G-Connected-M(G)}
Let $G$ be a disconnected graph. Then $M(G)$ is a connected graph if and only if $G$ has two components which are each either a $K_3$ or star graph.
\end{lem}

\begin{proof}
By direct computation we easily see that the matching complexes of  $K_3 \sqcup K_3$, $K_3 \sqcup S_n$, and $S_n \sqcup S_m$ are all connected graphs.

Suppose that $M(G)$ is a connected graph and $G$ is disconnected. Observe that if $G$ had more than two components or if any component contained two non-adjacent edges, then $G$ would contain a $3$-matching. Thus $G$ has exactly two components and all edges of each component are adjacent to each other, so the only possibility for each component of $G$ is $K_3$ or a star graph.
\end{proof}

We now will focus on \emph{connected} graphs whose matching complexes are also connected graphs. We start with the following lemma.

\begin{lem}\label{lem:path-4}
    Suppose $G$ is a connected graph with at least two edges. If $M(G)$ is also a connected graph, then $G$ contains a path of length four and no paths of length five or more.
\end{lem}

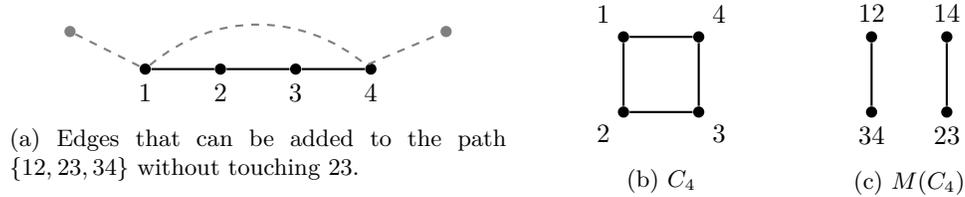
\begin{figure}
    \centering
    \begin{subfigure}{0.4\textwidth}
        \centering
        \begin{tikzpicture}
        \node (a) [vertex, label=below:$1$] at (0,0) {};
        \node (b) [vertex, label=below:$2$] at (1,0) {};
        \node (c) [vertex, label=below:$3$] at (2,0) {};
        \node (d) [vertex, label=below:$4$] at (3,0) {};
        
        \node (sa) [vertex, optional] at (-1,0.5) {};
        \node (sd) [vertex, optional] at (4,0.5) {};
        
        \draw [edge] (a) -- (b) -- (c) -- (d);
        
        \draw [edge, optional] (sa) -- (a) to [bend left=40] (d) -- (sd);
        \end{tikzpicture}
        \caption{Edges that can be added to the path $\set{12,23,34}$ without touching $23$.}\label{fig:short_path}
    \end{subfigure} \qquad
    \begin{subfigure}{0.15\textwidth}
        \centering
        \begin{tikzpicture}
        \node (a) [vertex, label=above left:$1$] at (0,1) {};
        \node (b) [vertex, label=below left:$2$] at (0,0) {};
        \node (c) [vertex, label=below right:$3$] at (1,0) {};
        \node (d) [vertex, label=above right:$4$] at (1,1) {};
        
        \draw [edge] (a) -- (b) -- (c) -- (d) -- (a);
        \end{tikzpicture}
        \caption{$C_4$}\label{fig:C_4}
    \end{subfigure} \qquad
    \begin{subfigure}{0.15\textwidth}
        \centering
        \begin{tikzpicture}
        \node (ab) [vertex, label=above:$12$] at (0,1) {};
        \node (bc) [vertex, label=below:$23$] at (1,0) {};
        \node (ad) [vertex, label=above:$14$] at (1,1) {};
        \node (cd) [vertex, label=below:$34$] at (0,0) {};
        
        \draw [edge] (ab) -- (cd) (ad) -- (bc);
        \end{tikzpicture}
        \caption{$M(C_4)$}\label{fig:M(C_4)}
    \end{subfigure}
    \caption{Graphs appearing in the proof of \cref{lem:path-4}}
    \label{fig:path-4}
\end{figure}

\begin{proof}
If $G$ has a path of length five or more, we get a 3-matching by taking the first, third, and fifth edges of the path, so $M(G)$ is not $1$-dimensional.

Assume that edges $e = uv$ and $e' = u'v'$ form a matching in $G$. Since $G$ is connected, there must be a path connecting $u$ to $u'$. Let $P$ be the shortest such path (which may contain $e$ or $e'$). Since the vertices of $e$ and $e'$ are distinct, $P$ must contain at least one edge besides $e$ and $e'$, so $P\cup\set{e,e'}$ is a path containing at least three edges. 

We now only need to show that $G$ has a path of length four in particular. Suppose $G$ contains a path of length three, say $\{12,23,34\}$. Now consider the edge $23$ (see \cref{fig:short_path}). Since $M(G)$ is a connected graph, there must be some other edge in $G$ which does not share endpoints with $23$, otherwise $23$ becomes an isolated vertex in $M(G)$.  Because $G$ is connected, the only way to do this without having a path of length four is to add the edge $14$, i.e., to have $C_4$ as a subgraph of our graph $G$. However, $M(C_4)$ is not a connected graph (see \cref{fig:C_4,fig:M(C_4)}), so this subgraph cannot contain all the edges of $G$. Furthermore, the matching complex of any graph on four vertices that contains a $4$-cycle is also disconnected. Thus we must have some edge in $G$ which (without loss of generality) shares no endpoints with both $12$ and $23$ while keeping $G$ a connected graph, which gives a path of length four.
\end{proof}

We now turn our attention again to \cref{fig:NP5k}, and note that each graph depicted in this figure contains a path of length four.

\begin{lem}\label{lem:connected-graph}
Suppose $G$ is a connected graph with at least two edges. Then $M(G)$ is a connected graph if and only if $G \in \mathcal{G}_1$, $\mathcal{G}_2$, or $\mathcal{G}_3$, or $G$ is the bowtie graph $B$.
\end{lem}

\begin{figure}
    \centering
    \begin{tikzpicture}[scale=1]
            \node (a) [vertex, label=left:$1$] at (0,0) {};
            \node (b) [vertex, label=above:$2$] at (1,0) {};
            \node (c) [vertex, label=below:$3$] at (1.5,-1) {};
            \node (d) [vertex, label=above:$4$] at (2,0) {};
            \node (e) [vertex, label=right:$5$] at (3,0) {};
            \node (sb1) [vertex, optional] at (0,0.7) {};
            \node (sd1) [vertex, optional] at (3,0.7) {};
            
            \draw [edge] (a)--(b)--(c)--(d)--(e);
            \draw [edge, highlight] (b)--(d);
            \draw [edge, optional] (b)--(sb1) (d)--(sd1);
        \end{tikzpicture}
    \caption{Graph appearing in Case~2(b) in the proof of \cref{lem:connected-graph}.}
    \label{fig:connected-graph}
\end{figure}
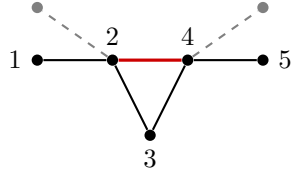

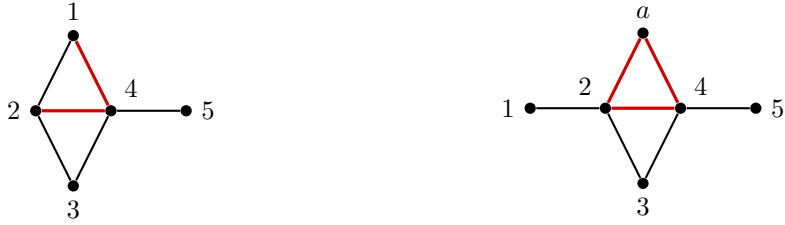
\begin{figure}
    \centering
    \begin{subfigure}{0.4\textwidth}
        \centering
        \begin{tikzpicture}[scale=1]
            \node (a) [vertex, label=above:$1$] at (1.5,1) {};
            \node (b) [vertex, label=left:$2$] at (1,0) {};
            \node (c) [vertex, label=below:$3$] at (1.5,-1) {};
            \node (d) [vertex, label=above right:$4$] at (2,0) {};
            \node (e) [vertex, label=right:$5$] at (3,0) {};
            
            \draw [edge] (a)--(b)--(c)--(d)--(e);
            \draw [edge, highlight] (b)--(d) (a)--(d);
        \end{tikzpicture}
    \end{subfigure} \quad
    \begin{subfigure}{0.4\textwidth}
        \centering
        \begin{tikzpicture}[scale=1]
            \node (a) [vertex, label=left:$1$] at (0,0) {};
            \node (b) [vertex, label=above left:$2$] at (1,0) {};
            \node (c) [vertex, label=below:$3$] at (1.5,-1) {};
            \node (d) [vertex, label=above right:$4$] at (2,0) {};
            \node (e) [vertex, label=right:$5$] at (3,0) {};
            \node (sbd1) [vertex, label=above:$a$] at (1.5,1) {};
            
            \draw [edge] (a)--(b)--(c)--(d)--(e);
            \draw [edge, highlight] (b)--(d) (b)--(sbd1)--(d);
        \end{tikzpicture}
    \end{subfigure}
    \caption{Graphs appearing in Case~4(b) in the proof of \cref{lem:connected-graph}.}
    \label{fig:Case-4(b)}
\end{figure}

\begin{proof}
It is straightforward to check that the matching complex of any graph in these families is indeed a connected graph; we omit the details of these calculations.

For the other direction of the proof, assume $M(G)$ is a connected graph. By \cref{lem:path-4}, $G$ contains a path $P$ of length four, say $P = \{12,23,34,45\}$. If $G$ is simply this path, then $G \in \mathcal{G}_1$. If this path is not all of $G$, let us consider what we can add. Note that we cannot add any edges which produce a 3-matching in the graph, so in particular we cannot have any edges of the form $1a$, $3a$, or $5a$ where $a$ is a new vertex, and we cannot add any edges that do not share a vertex with this path of length four. We will break the remainder of the proof into cases.

\textbf{Case~1:} Suppose $G$ has $C_5$ as a subgraph. We note that $M(C_5)$ is itself a connected graph (namely, $M(C_5) = C_5$). If $G$ contains any edge whose endpoints are not both contained in this $C_5$, then $G$ contains a 3-matching. However, $M(G)$ remains connected and $1$-dimensional if we add any number of edges between vertices in the $C_5$. Thus $G$ is in the family $\mathcal{G}_3$.

\textbf{Case 2:} Suppose $G$ has $C_3$ as a subgraph but not $C_4$ or $C_5$. Without loss of generality, this can only occur if we add edge $13$ or $24$ to our path of length four because otherwise we immediately get a 3-matching. We will consider these as Cases~2(a) and 2(b), respectively.

In Case~2(a), we add the edge $13$. Adding any edge to vertices $1$ or $2$ will either create a $4$-cycle (if no new vertices are introduced) or $3$-matching (if there are new vertices) and thus isn't allowed. Adding any edge to vertex $5$ except the edge $35$ will again create a $4$-cycle or $3$-matching. Thus the only allowed options are to add any number of pendant edges off vertex $4$ or to instead add the edge $35$. Observe that doing both of these would create a $3$-matching. The first of these options puts $G$ in $\mathcal{G}_2$, the second shows that $G$ is the bowtie graph $B$.

In Case~2(b), we add the edge $24$ as depicted in \cref{fig:connected-graph}. Observe that any additional edge with $1$, $3$, or $5$ as an endpoint will create either $4$-cycle or $3$-matching and thus is not allowed. The only possible additional edges in this case have either $2$ or $4$ as an endpoint. However, any graph in this family has a disconnected matching complex---in particular, the edge $24$ is not in a matching with any other edge. Thus this case is impossible.

\textbf{Case 3:} Suppose $G$ has $C_4$ as a subgraph but not $C_3$ or $C_5$. Without loss of generality, the only way this can occur without introducing a 3-matching is to add edge $14$ to the $P_4$ subgraph. Observe that any additional edge with $1$ or $3$ as an endpoint will either create a disallowed cycle or $3$-matching. This is the same for any edge with $5$ as an endpoint except the edge $25$, which is allowed. The edge $24$ would create a $3$-cycle and thus is not allowed. However, any number of pendants off $2$ and $4$ and any paths of length two connecting vertices $2$ and $4$ are allowed. Therefore $G \in \mathcal{G}_1$.

\textbf{Case 4:} Suppose $G$ has both $C_3$ and $C_4$ as a subgraph but not $C_5$. There are two ways to introduce the $C_3$ (without loss of generality): As in Case~2, we can add the edge $13$ or $24$. We will call these Cases~4(a) and 4(b), respectively.

In Case~4(a), we add the edge $13$. The only edges involving $1$ and $2$ that we can add without introducing a $3$-matching or $5$-cycle are $14$ and $24$. We must add at least one of these to create a $4$-cycle in $G$. Once we do so, we may add the other edge and any number of pendants off $4$. Any edge with $5$ as an endpoint will create a $3$-matching or $5$-cycle. Thus $G \in \mathcal{G}_2$.

In Case~4(b), we add the edge $24$. As before, pendant edges off $1$, $3$, or $5$ create $3$-matchings, and adding the edge $15$ creates a $5$-cycle. Without loss of generality, the only way to create a $4$-cycle without introducing new vertices is to add the edge $14$. The only way to prevent $24$ from being an isolated vertex in the matching complex is to add the edge $13$. Now adding pendants to $2$ produces a $3$ matching, thus $G \in \mathcal{G}_2$.

If instead we create a $4$-cycle with a new vertex, the only possible way without creating $3$-matchings or disallowed cycles is to add $2a$ and $4a$ for a new vertex $a$ as in \cref{fig:Case-4(b)}. However, as in Case~2(b), we see that $24$ must be an isolated vertex in the matching complex. Thus this case is impossible.

\textbf{Case 5:} Suppose $G$ contains no cycles. In this case the only edges we can add without getting a 3-matching are pendants off vertices $2$ and $4$. Thus $G \in \mathcal{G}_1$.
\end{proof}

We immediately get the following corollary, combining the above results of this section with the definition of a two-dimensional Buchsbaum complex.

\begin{cor}\label{cor:Easy-Buchs}
Let $G$ be a graph. Then $M(G)$ is a two-dimensional Buchsbaum complex if and only if for each edge $e$ of $G$, either 
\begin{enumerate}[(a)]
    \item $N_e$ has two components which are each either a $K_3$ or a star graph, or
    \item $N_e$ is in one of the families $\mathcal{G}_1$, $\mathcal{G}_2$, and $\mathcal{G}_3$, or $N_e$ is the bowtie graph $B$.
    \end{enumerate}
(That is, $N_e$ is one of the graphs described in \cref{prop:CM-Buchs-1-dim}(a).)
\end{cor}

We are now able to complete the proof of \cref{prop:CM-Buchs-1-dim}, characterizing graphs whose matching complex is one-dimensional and either Cohen--Macaulay or Buchsbaum.

\begin{proof}[Proof of \cref{prop:CM-Buchs-1-dim}]
We recall that a $1$-dimensional complex is Cohen--Macaulay if and only if it is connected, and it is Buchsbaum if and only if it has no isolated vertices.

Therefore \cref{prop:CM-Buchs-1-dim}(a) follows immediately. For \cref{prop:CM-Buchs-1-dim}(b), the only graphs that have disconnected matching complexes with no isolated vertices are $K_4$ and $C_4$ by \cite[Theorem 2.9]{BGM}.
\end{proof}

\section{Buchsbaum graph families}\label{sec:Families}

The goal of this section is to provide an explicit description of all graphs $G$ for which $M(G)$ is two-dimensional and Buchsbaum. The following is our main result.

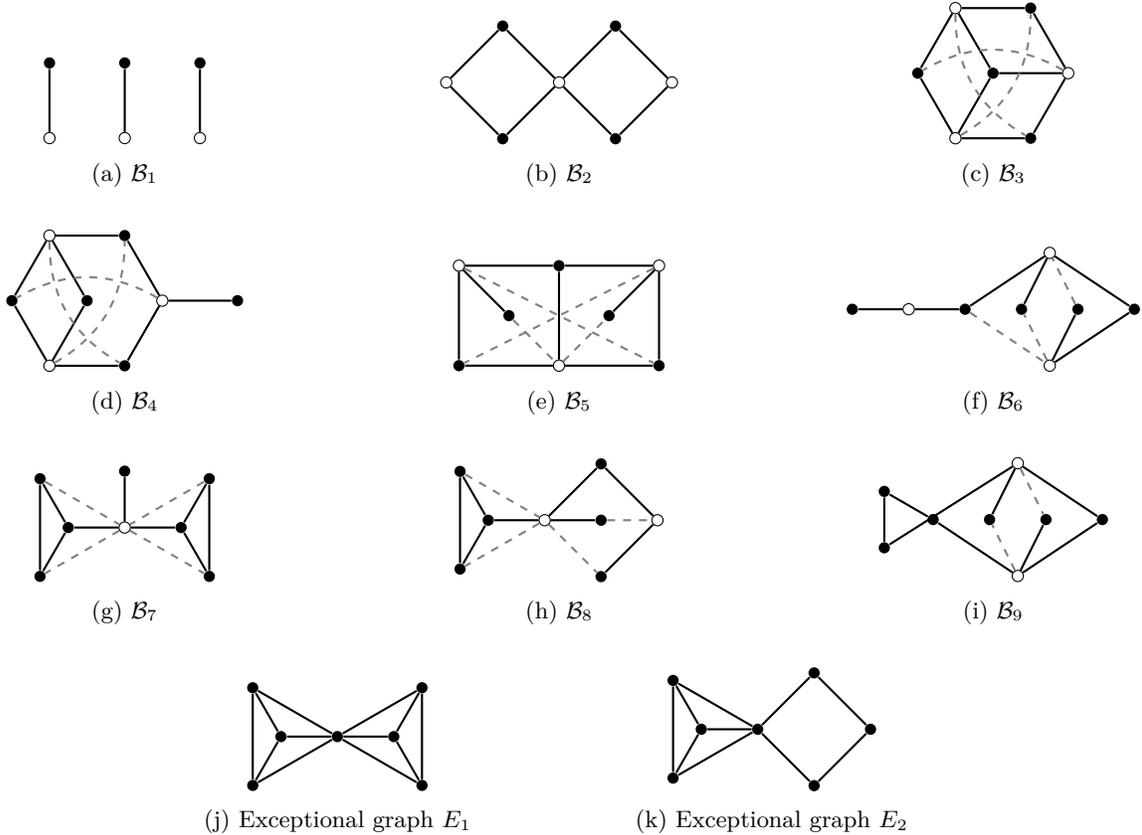
\begin{figure}
    \centering
    \begin{subfigure}[b]{0.3\textwidth}
        \centering
        \begin{tikzpicture}
        \node (g1) [vertex, starpoint] at (0,0) {};
        \node (g2) [vertex, starpoint] at (1,0) {};
        \node (g3) [vertex, starpoint] at (2,0) {};
        \node (v1) [vertex] at (0,1) {};
        \node (v2) [vertex] at (1,1) {};
        \node (v3) [vertex] at (2,1) {};
        \draw [edge] (g1)--(v1) (g2)--(v2) (g3)--(v3);
        \end{tikzpicture}
        \caption{$\mathcal{B}_1$}
    \end{subfigure} \qquad
    \begin{subfigure}[b]{0.3\textwidth}
        \centering
        \begin{tikzpicture}[scale=0.75]
        \node (g1) [vertex, starpoint] at (0,0) {};
        \node (g2) [vertex, starpoint] at (2,0) {};
        \node (g3) [vertex, starpoint] at (4,0) {};
        \node (v1) [vertex] at (1,1) {};
        \node (v2) [vertex] at (3,1) {};
        \node (v3) [vertex] at (1,-1) {};
        \node (v4) [vertex] at (3,-1) {};
        \draw [edge] (g1)--(v1)--(g2)--(v2)--(g3)--(v4)--(g2)--(v3)--(g1);
        \end{tikzpicture}
        \caption{$\mathcal{B}_2$}
    \end{subfigure} \qquad
    \begin{subfigure}[b]{0.3\textwidth}
        \centering
        \begin{tikzpicture}[scale=1]
        \node (v1) [vertex] at (0:0) {};
        \node (g1) [vertex, starpoint] at (0:1) {};
        \node (v2) [vertex] at (-60:1) {};
        \node (g2) [vertex, starpoint] at (-120:1) {};
        \node (v3) [vertex] at (-180:1) {};
        \node (g3) [vertex, starpoint] at (-240:1) {};
        \node (v4) [vertex] at (-300:1) {};
        \draw [edge] (g1)--(v2)--(g2)--(v3)--(g3)--(v4)--(g1) (v1)--(g1) (v1)--(g2) (v1)--(g3);
        \draw [edge, optional] (g1) to [bend right] (v3) (g2) to [bend right] (v4) (g3) to [bend right] (v2);
        \end{tikzpicture}
        \caption{$\mathcal{B}_3$}
    \end{subfigure} \\[0.5cm]
    \begin{subfigure}[b]{0.3\textwidth}
        \centering
        \begin{tikzpicture}[scale=1]
        \node (v1) [vertex] at (0:0) {};
        \node (g1) [vertex, starpoint] at (0:1) {};
        \node (v2) [vertex] at (-60:1) {};
        \node (g2) [vertex, starpoint] at (-120:1) {};
        \node (v3) [vertex] at (-180:1) {};
        \node (g3) [vertex, starpoint] at (-240:1) {};
        \node (v4) [vertex] at (-300:1) {};
        \node (v5) [vertex] at (0:2) {};
        \draw [edge] (g1)--(v2)--(g2)--(v3)--(g3)--(v4)--(g1) (g1)--(v5) (g2)--(v1)--(g3);
        \draw [edge, optional] (g1) to [bend right] (v3) (g2) to [bend right] (v4) (g3) to [bend right] (v2);
        \end{tikzpicture}
        \caption{$\mathcal{B}_4$}
    \end{subfigure} \qquad
    \begin{subfigure}[b]{0.3\textwidth}
        \centering
        \begin{tikzpicture}[scale=1.33, rotate=-90]
        \node (g1) [vertex, starpoint] at (0,1) {};
        \node (v1) [vertex] at (1,1) {};
        \node (v2) [vertex] at (0.5,0.5) {};
        \node (v3) [vertex] at (0,0) {};
        \node (g2) [vertex, starpoint] at (1,0) {};
        \node (v4) [vertex] at (0.5,-0.5) {};
        \node (g3) [vertex, starpoint] at (0,-1) {};
        \node (v5) [vertex] at (1,-1) {};
        \draw [edge] (v2)--(g1)--(v1)--(g2)--(v3)--(g1) (v3)--(g3)--(v4) (g3)--(v5)--(g2);
        \draw [edge, optional] (v4)--(g2)--(v2) (g1) -- (v5) (g3) -- (v1);
        \end{tikzpicture}
        \caption{$\mathcal{B}_5$}
    \end{subfigure} \qquad
    \begin{subfigure}[b]{0.3\textwidth}
        \centering
        \begin{tikzpicture}[scale=0.75]
        \node (v1) [vertex] at (0,0) {};
        \node (g1) [vertex, starpoint] at (1,0) {};
        \node (v2) [vertex] at (2,0) {};
        \node (v3) [vertex] at (3,0) {};
        \node (v4) [vertex] at (4,0) {};
        \node (v5) [vertex] at (5,0) {};
        \node (g2) [vertex, starpoint] at (3.5,1) {};
        \node (g3) [vertex, starpoint] at (3.5,-1) {};
        \draw [edge] (v1)--(g1)--(v2)--(g2)--(v5)--(g3)--(v4) (g2)--(v3);
        \draw [edge, optional] (v2)--(g3)--(v3) (v4)--(g2);
        \end{tikzpicture}
        \caption{$\mathcal{B}_6$}
    \end{subfigure} \\[0.5cm]
    \begin{subfigure}[b]{0.3\textwidth}
        \centering
        \begin{tikzpicture}[scale=0.75]
        \node (g) [vertex, starpoint] at (0,0) {};
        \node (v1) [vertex] at (-1,0) {};
        \node (v2) [vertex] at ($(v1) + (120:1)$) {};
        \node (v3) [vertex] at ($(v1) + (240:1)$) {};
        \node (v4) [vertex] at (1,0) {};
        \node (v5) [vertex] at ($(v4) + (60:1)$) {};
        \node (v6) [vertex] at ($(v4) + (-60:1)$) {};
        \node (v7) [vertex] at (0,1) {};
        \draw [edge] (v1)--(v2)--(v3)--(v1)--(g)--(v4)--(v5)--(v6)--(v4) (g)--(v7);
        \draw [edge, optional] (v2)--(g)--(v3) (v5)--(g)--(v6);
        \end{tikzpicture}
        \caption{$\mathcal{B}_7$}
    \end{subfigure} \qquad
    \begin{subfigure}[b]{0.3\textwidth}
        \centering
        \begin{tikzpicture}[scale=0.75]
        \node (g1) [vertex, starpoint] at (0,0) {};
        \node (v1) [vertex] at (-1,0) {};
        \node (v2) [vertex] at ($(v1) + (120:1)$) {};
        \node (v3) [vertex] at ($(v1) + (240:1)$) {};
        \node (v4) [vertex] at (1,1) {};
        \node (v5) [vertex] at (1,0) {};
        \node (v6) [vertex] at (1,-1) {};
        \node (g2) [vertex, starpoint] at (2,0) {};
        \draw [edge] (v1)--(v2)--(v3)--(v1)--(g1)--(v4)--(g2)--(v6) (g1)--(v5);
        \draw [edge, optional] (v2)--(g1)--(v3) (v5)--(g2) (g1)--(v6);
        \end{tikzpicture}
        \caption{$\mathcal{B}_8$}
    \end{subfigure} \qquad
    \begin{subfigure}[b]{0.3\textwidth}
        \centering
        \begin{tikzpicture}[scale=0.75]
        \node (v1) [vertex] at (0,0) {};
        \node (v2) [vertex] at (150:1) {};
        \node (v3) [vertex] at (210:1) {};
        \node (v4) [vertex] at (1,0) {};
        \node (v5) [vertex] at (2,0) {};
        \node (v6) [vertex] at (3,0) {};
        \node (g1) [vertex, starpoint] at (1.5,1) {};
        \node (g2) [vertex, starpoint] at (1.5,-1) {};
        \draw [edge] (g2)--(v1)--(v2)--(v3)--(g1)--(v6)--(g2)--(v5) (g1)--(v4);
        \draw [edge, optional] (g2)--(v4) (g1)--(v5);
        \end{tikzpicture}
        \caption{$\mathcal{B}_9$}
    \end{subfigure} \\[0.5cm]
    \begin{subfigure}[b]{0.3\textwidth}
        \centering
        \begin{tikzpicture}[scale=0.75]
        \node (g) [vertex] at (0,0) {};
        \node (v1) [vertex] at (-1,0) {};
        \node (v2) [vertex] at ($(v1) + (120:1)$) {};
        \node (v3) [vertex] at ($(v1) + (240:1)$) {};
        \node (v4) [vertex] at (1,0) {};
        \node (v5) [vertex] at ($(v4) + (60:1)$) {};
        \node (v6) [vertex] at ($(v4) + (-60:1)$) {};
        \draw [edge] (v1)--(v2)--(v3)--(v1)--(g)--(v4)--(v5)--(v6)--(v4);
        \draw [edge] (v2)--(g)--(v3) (v5)--(g)--(v6);
        \end{tikzpicture}
        \caption{Exceptional graph $E_1$}
    \end{subfigure} \qquad
    \begin{subfigure}[b]{0.3\textwidth}
        \centering
        \begin{tikzpicture}[scale=0.75]
        \node (g1) [vertex] at (0,0) {};
        \node (v1) [vertex] at (-1,0) {};
        \node (v2) [vertex] at ($(v1) + (120:1)$) {};
        \node (v3) [vertex] at ($(v1) + (240:1)$) {};
        \node (v4) [vertex] at (1,1) {};
        \node (v6) [vertex] at (1,-1) {};
        \node (g2) [vertex] at (2,0) {};
        \draw [edge] (v1)--(v2)--(v3)--(v1)--(g1)--(v4)--(g2)--(v6);
        \draw [edge] (v2)--(g1)--(v3) (g1)--(v6);
        \end{tikzpicture}
        \caption{Exceptional graph $E_2$}
    \end{subfigure} \qquad
    \caption{Graphs whose matching complexes are two-dimensional and Buchsbaum. As before, solid edges are required and dotted edges are optional. For each graph, any number of filled vertices may be introduced and each new vertex attached to some nonzero number of unfilled vertices.
    \label{fig:2D-buchsbaum-graphs}}
\end{figure}

\begin{thm}\label{thm:MainTheorem}
Let $G$ be a graph. Then $M(G)$ is a two-dimensional Buchsbaum complex if and only if $G$ is one of the following graphs (which are defined below and depicted in \cref{fig:2D-buchsbaum-graphs}):
\begin{enumerate}[(a)]
\item a member of one of the families $\mathcal{B}_{C_7},$ $\mathcal{B}_P$, or $\mathcal{B}_i$ for some $i \in \{1, \dotsc, 9\}$,
\item one of the two exceptional graphs $E_1$ and $E_2$, or
\item one of the disconnected graphs described in \cref{prop:Disconnected2dBuchs}, i.e.,
\begin{enumerate}[(i)]
    \item $G$ has three components, each of which is either $K_3$ or a star graph, or
    \item $G$ has two components, one of which is $K_3$ or a star graph and the other is either the bowtie graph $B$ or a graph in one of the families $\mathcal{G}_1$, $\mathcal{G}_2$, and $\mathcal{G}_3$.
\end{enumerate}
\end{enumerate}
\end{thm}

The rest of this section is devoted to proving \cref{thm:MainTheorem}, and we will now provide a brief outline of the proof. We first briefly consider when $G$ is disconnected in \cref{prop:Disconnected2dBuchs}. Next, we collect a variety of results on graphs containing cycles of certain sizes, including bipartite graphs. We then split up our remaining casework using the notion of a ``link connected'' graph, which we introduce in \cref{def:link_conn} (a graph is link connected if all non-adjacent subgraphs $N_e$ are connected). In \cref{sec:Link-Connected}, we examine the graphs that are link connected: We consider which cycles can appear in these graphs, and use this analysis to deduce which families these graphs belong to.  
This leaves the non-link-connected graphs, which we deal with in \cref{sec:Not-Link-Connected}. These graphs by definition have some non-adjacent subgraph $N_e$ that is not connected. By \cref{cor:Easy-Buchs}, there are only a few possibilities for what this $N_e$ can be, and we examine each possibility one by one.

First, let us handle the case where $G$ is disconnected.

\begin{prop}\label{prop:Disconnected2dBuchs}
Suppose $G$ is a disconnected graph. Then $M(G)$ is two-dimensional and Buchsbaum if and only if either
\begin{enumerate}[(a)]
    \item $G$ has three components, each of which is either $K_3$ or a star graph, or \label{prop:Disconnected2dBuchs-3comp}
    \item $G$ has two components, one of which is $K_3$ or a star graph and the other is either the bowtie graph $B$ or a graph in one of the families $\mathcal{G}_1$, $\mathcal{G}_2$, and $\mathcal{G}_3$. \label{prop:Disconnected2dBuchs-2comp}
\end{enumerate}
\end{prop}

\begin{proof}
First, if $G$ is a graph satisfying \ref{prop:Disconnected2dBuchs-3comp} or \ref{prop:Disconnected2dBuchs-2comp}, it is straightforward to check that $M(G)$ is indeed two-dimensional and Buchsbaum.

Conversely, suppose that $G$ is a disconnected graph and $M(G)$ is two-dimensional and Buchsbaum. If $G$ has four or more connected components, then it is guaranteed to have a 4-matching and $M(G)$ is thus not two-dimensional. If $G$ has three components, then each must be either $K_3$ or a star graph because each component must have all edges adjacent to each other; otherwise we could find a $4$-matching in $G$.

Next, suppose that $G$ has exactly two connected components and $G = G_1 \sqcup G_2$. At least one component must not contain a $2$-matching: Say this component is $G_1$. Then $G_1$ is $K_3$ or a star graph. If we take any edge $e \in G_1$, we see that $N_e$ is precisely $G_2$. Thus by \cref{lem:connected-graph}, $G_2$ must be in $\mathcal G_1$, $\mathcal G_2$, $\mathcal G_3$ or $G_2$ is the bowtie graph $B$.
\end{proof}

We now turn our attention to \cref{fig:2D-buchsbaum-graphs}, which depicts several families of graphs $\mathcal{B}_i$ and two exceptional graphs $E_1$ and $E_2$ which all have two-dimensional Buchsbaum matching complexes. The families are defined as follows: Solid edges are necessary and dotted edges are optional. Observe the filled and unfilled vertices in \cref{fig:2D-buchsbaum-graphs}: For each graph, we may add any number of new filled vertices and attach each to a non-zero number of the unfilled vertices. We note that the families $\mathcal{B}_1$ through $\mathcal{B}_6$ have nonempty intersection. Otherwise these families are mutually disjoint.

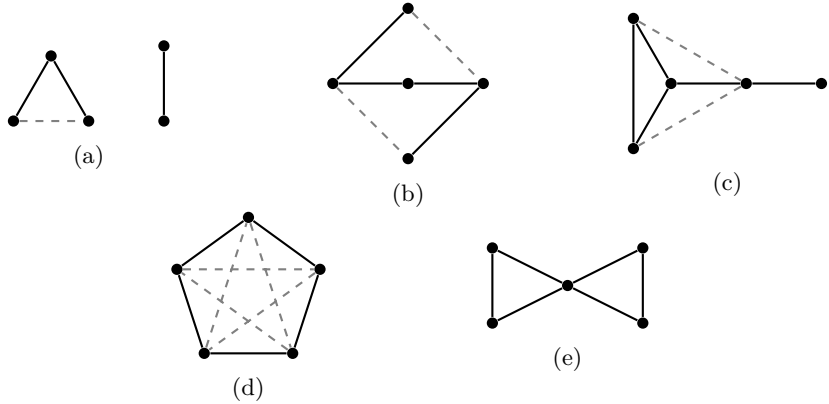
\begin{figure}
    \centering
    \begin{subfigure}{0.25\textwidth}
        \centering
        \begin{tikzpicture}
            \node (a) [vertex] at (0,0) {};
            \node (b) [vertex] at (60:1) {};
            \node (c) [vertex] at (0:1) {};
            \node (d) [vertex] at (2,0) {};
            \node (e) [vertex] at (2,1) {};
            
            \draw [edge] (a)--(b)--(c) (d)--(e);
            \draw [edge, optional] (a)--(c);
        \end{tikzpicture}
        \caption{}
    \end{subfigure}
    \begin{subfigure}{0.25\textwidth}
        \centering
        \begin{tikzpicture}
            \node (a) [vertex] at (1,1) {};
            \node (b) [vertex] at (0,0) {};
            \node (c) [vertex] at (1,0) {};
            \node (d) [vertex] at (2,0) {};
            \node (e) [vertex] at (1,-1) {};
            
            \draw [edge] (a)--(b)--(c)--(d)--(e);
            \draw [edge, optional] (a)--(d) (b)--(e);
        \end{tikzpicture}
        \caption{}
    \end{subfigure}
    \begin{subfigure}{0.25\textwidth}
        \centering
        \begin{tikzpicture}[scale=1]
            \node (a) [vertex] at (120:1) {};
            \node (b) [vertex] at (240:1) {};
            \node (c) [vertex] at (0,0) {};
            \node (d) [vertex] at (0:1) {};
            \node (e) [vertex] at (0:2) {};
            
            \draw [edge] (a)--(b)--(c)--(d)--(e) (a)--(c);
            \draw [edge, optional] (a)--(d)--(b);
        \end{tikzpicture}
        \caption{}
    \end{subfigure}
    \begin{subfigure}{0.25\textwidth}
        \centering
        \begin{tikzpicture}[scale=1]
            \node (a) [vertex] at (90:1) {};
            \node (b) [vertex] at (18:1) {};
            \node (c) [vertex] at (306:1) {};
            \node (d) [vertex] at (234:1) {};
            \node (e) [vertex] at (162:1) {};
            
            \draw [edge] (a)--(b)--(c)--(d)--(e)--(a);
            \draw [edge, optional] (a)--(c)--(e)--(b)--(d)--(a);
        \end{tikzpicture}
        \caption{}
    \end{subfigure}
    \begin{subfigure}{0.25\textwidth}
        \centering
        \begin{tikzpicture}[scale=1]
            \node (a) [vertex] at (0,1) {};
            \node (b) [vertex] at (0,0) {};
            \node (c) [vertex] at (1,0.5) {};
            \node (d) [vertex] at (2,1) {};
            \node (e) [vertex] at (2,0) {};
            
            \draw [edge] (a)--(b)--(c)--(d)--(e);
            \draw [edge] (a)--(c)--(e);
        \end{tikzpicture}
        \caption{}
    \end{subfigure}
    \caption{Graphs with $5$ vertices whose matching complex is a connected graph.}
    \label{fig:5-verts-M(G)-connected}
\end{figure}

We now describe two additional families that are not depicted in \cref{fig:2D-buchsbaum-graphs}. The first is $\mathcal{B}_{C_7}$, which is defined to be the family of all graphs containing $C_7$ as a subgraph that have two-dimensional Buchsbaum matching complexes. We discuss this family in more detail later in this section, in particular in \cref{lem:c7-no-extra-vertices,lem:C7links,tab:c7-data,fig:5-verts-M(G)-connected}.

\begin{figure}
        \centering
        \begin{tikzpicture}[rotate=17, scale=0.7]
        \node (c) [vertex] at (0,0) {};
        \node (s1c) [vertex] at ($(c) + (0:1.6)$) {};
        \foreach \i in {1,...,6}
            \node (s1\i) [vertex] at ($(s1c) + (180+\i*360/7:1)$) {};
        \node (s2c) [vertex] at ($(c) + (120:1.6)$) {};
        \foreach \i in {1,...,4}
            \node (s2\i) [vertex] at ($(s2c) + (300+\i*360/5:1)$) {};
        \node (t1) [vertex] at ($(c) + (240+25:1.8)$) {};
        \node (t2) [vertex] at ($(c) + (240-25:1.8)$) {};
        \draw [edge] (s1c)--(c)--(t1)--(t2)--(c)--(s2c) foreach \i in {1,...,6} {(s1c)--(s1\i)} foreach \j in {1,...,4} {(s2c)--(s2\j)};
        \end{tikzpicture}
        \caption{Example of a petal graph in $\mathcal{B}_P$. Each petal is either a $K_3$ or star graph with at least two edges.} \label{fig:petal}
    \end{figure}
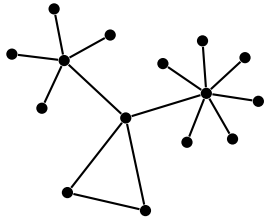

The last family is $\mathcal{B}_P$, which we call petal graphs. These are formed by taking three graphs---each either a $K_3$ or star graph with at least two edges---and then gluing these graphs together at a single vertex. For the star graphs, the gluing vertex must be a non-central vertex of the star. The resulting graph will have one main central vertex and three `petals,' each of which is a $K_3$ or star graph. For an example of a petal graph, see \cref{fig:petal}.

We can straightforwardly verify that any graph in \cref{fig:2D-buchsbaum-graphs} or the families $\mathcal{B}_{C_7}$ and $\mathcal{B}_P$ has a two-dimensional Buchsbaum matching complex: For petal graphs and the families in \cref{fig:2D-buchsbaum-graphs}, we observe that---regardless of whether any dotted edges are included---the non-adjacent subgraph of any edge is a graph from \cref{prop:CM-Buchs-1-dim}. Thus we only need to consider added edges in \cref{fig:2D-buchsbaum-graphs}. The non-adjacent subgraphs of these edges correspond precisely to removing an unfilled vertex from the graph in question; again, we see that all such graphs appear in \cref{prop:CM-Buchs-1-dim}.

We will spend the remainder of this section showing that any connected graph with a two-dimensional Buchsbaum matching complex must be in one of the families in \cref{fig:2D-buchsbaum-graphs} or the $\mathcal{B}_{C_7}$ or $\mathcal{B}_P$ families. Doing so will complete the proof of \cref{thm:MainTheorem}.

We recall that the only bipartite graphs in \cref{lem:Disconnected-G-Connected-M(G),lem:connected-graph} are either the disjoint union of two star graphs or in graph family $\mathcal{G}_1$. Considering \cref{fig:2D-buchsbaum-graphs}, we make precise the following observation.

\begin{prop}\label{prop:Rowan}
If $G$ is a connected, bipartite graph and $M(G)$ is a two-dimensional Buchsbaum complex, then one side of the bipartition has exactly three vertices.
\end{prop}

\begin{proof}
First, since $M(G)$ is two-dimensional, there is some matching of three edges in $G$. Each of these edges must have one vertex in each side of the bipartition, so each side must contain at least three vertices.

Now, let $e = ab$ be any edge of $G$. We know that the matching complex of $N_e$ is a connected graph and that $N_e$ itself must be bipartite, so by \cref{lem:Disconnected-G-Connected-M(G),lem:connected-graph}, $N_e$ is either a graph in $\mathcal{G}_1$ or a disjoint union of two star graphs.

Assume $N_e \in \mathcal{G}_1$. For every graph in $\mathcal{G}_1$, any bipartition has one side with only two vertices, namely the two vertices labeled $2$ and $4$ in \cref{fig:NP5k}. The edge $e$ itself contributes one more vertex to each side, so this side has exactly $3$ vertices in $G$.

This only leaves the case where $N_e$ is a disjoint union of two stars, i.e., $N_e = S_m \sqcup S_n$. Since this subgraph is disconnected, there is more than one way to bipartition it---the centers of the stars can be on either the same side of the bipartition or opposite sides. If the centers of the stars are on the same side (or either of $i$ and $j$ is $1$), the same argument as above works, since $e$ again contributes one more vertex to this side.
However, in the bipartition of $N_e$ which puts the centers on opposite sides (see \cref{fig:Rowan}), each side may be arbitrarily large: We must argue that this is not allowed in $G$.

Since $G$ is connected, there must be at least one edge in $G$ connecting $e$ to each star. We claim that such an edge cannot connect to the center of a star: Assume that the edge $x$ connects vertex $a$ to the center of $S_m$. Given the constraints on $G$, we see that $N_x$ cannot contain a path of length four, so $N_x$ must be disconnected. Thus there cannot be an edge between $b$ and the center of $S_n$. Therefore we must have an edge between $a$ and a non-central vertex of $S_n$; call this edge $y$. Returning to $N_x$, we see that there must be an edge containing vertex $b$ that does not connect to $S_n$; call this edge $z$. Observe that this creates a $4$-matching in $G$: Take the edges $y$ and $z$, any edge in $S_m$ that is not adjacent to $z$, and any edge in $S_n$ that is not adjacent to $y$.

Similarly, if there are only edges between $e$ and non-central vertices of $S_m$ and $S_n$ then $G$ will also contain a $4$-matching: Take two such non-adjacent edges, then for each star we can always find another non-adjacent edge.

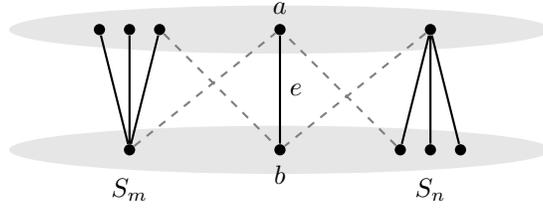
\begin{figure}
    \centering
    \begin{tikzpicture}[scale=0.8]
    \fill [black!10] (3,2) ellipse [x radius=4.5cm, y radius=0.4cm];
    \fill [black!10] (3,0) ellipse [x radius=4.5cm, y radius=0.4cm];
    \node (s1) [vertex] at (0,2) {};
    \node (s2) [vertex] at (0.5,2) {};
    \node (s3) [vertex] at (1,2) {};
    \node (sc) [vertex, label={[label distance=5pt]below:$S_m$}] at (0.5,0) {};
    \node (e1) [vertex, label=above:$a$] at (3,2) {};
    \node (e2) [vertex, label=below:$b$] at (3,0) {};
    \node (t1) [vertex] at (5,0) {};
    \node (t2) [vertex, label={[label distance=5pt]below:$S_n$}] at (5.5,0) {};
    \node (t3) [vertex] at (6,0) {};
    \node (tc) [vertex] at (5.5,2) {};
    \draw [edge] (s1)--(sc) (s2)--(sc) (s3)--(sc) (e1) -- node {$e$} (e2) (t1)--(tc) (t2)--(tc) (t3)--(tc);
    \draw [edge, optional] (s3)--(e2) (sc)--(e1) (t1)--(e1) (tc)--(e2);
    \end{tikzpicture}
    \caption{The case where $N_e$ is a disjoint union of two stars, bipartitioned badly, in \cref{prop:Rowan}.}
    \label{fig:Rowan}
\end{figure}

Thus this case is impossible, so one side of the bipartition of $G$ must have exactly $3$ vertices.
\end{proof}

There are several relevant observations about non-bipartite graphs that we can make as well.

\begin{prop}\label{prop:C_5-implies-C_7}
    Assume $G$ is connected and $M(G)$ is a two-dimensional Buchsbaum complex. If $G$ contains $C_5$, then $G$ contains $C_7$.
\end{prop}

\begin{proof}
Suppose first that $G$ contains a $5$-cycle and an edge $e=ab$ that is disjoint from this cycle, as depicted in \cref{fig:C_5PlusEdge}.

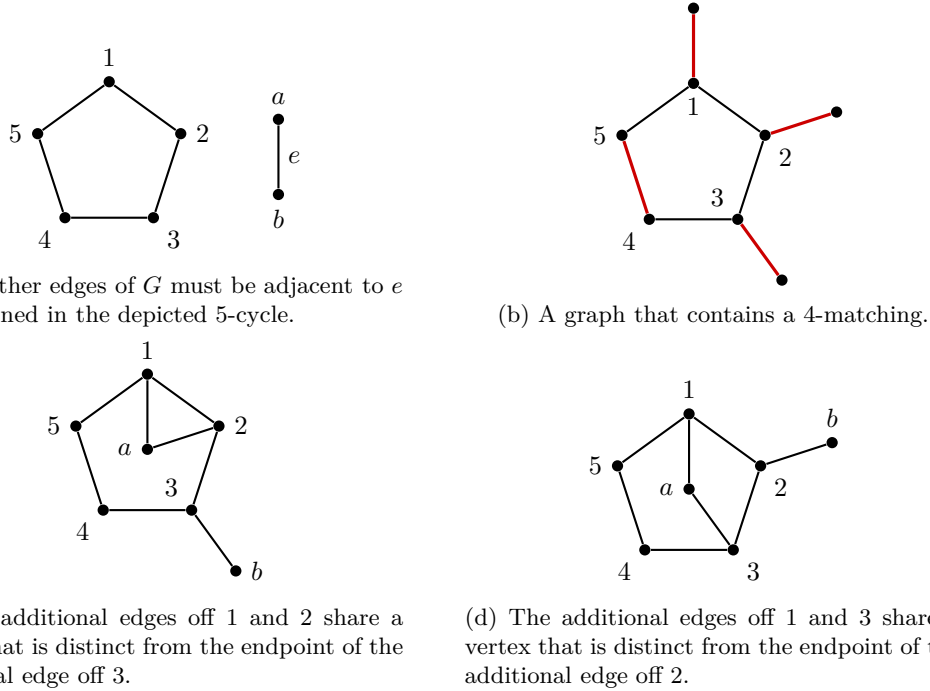
\begin{figure}
    \centering
    \begin{subfigure}[b]{0.4\textwidth}
        \centering
        \begin{tikzpicture}
        \pgftransformrotate{90}
        \node (4) [vertex, label=below left:$4$] at (144:1) {};
        \node (5) [vertex, label=left:$5$] at (72:1) {};
        \node (1) [vertex, label=above:$1$] at (0:1) {};
        \node (2) [vertex, label=right:$2$] at (288:1) {};
        \node (3) [vertex, label=below right:$3$] at (216:1) {};

        \node (a) [vertex, label=above:$a$] at (.5,-2.25) {};
        \node (b) [vertex, label=below:$b$] at (-.5,-2.25) {};

        \draw [edge] (1) -- (2) -- (3) -- (4) -- (5) -- (1);
        \draw [edge] (a) to node [auto] {$e$} (b);
        \end{tikzpicture}
        \caption{All other edges of $G$ must be adjacent to $e$ or contained in the depicted $5$-cycle.}
        \label{fig:C_5PlusEdge}
    \end{subfigure} \qquad
    \begin{subfigure}[b]{0.4\textwidth}
        \centering
        
        \begin{tikzpicture}
        \pgftransformrotate{90}
        \node (4) [vertex, label=below left:$4$] at (144:1) {};
        \node (5) [vertex, label=left:$5$] at (72:1) {};
        \node (1) [vertex, label=below:$1$] at (0:1) {};
        \node (2) [vertex, label=below right:$2$] at (288:1) {};
        \node (3) [vertex, label=above left:$3$] at (216:1) {};

        \node (a) [vertex] at (0:2) {};
        \node (b) [vertex] at (288:2) {};
        \node (c) [vertex] at (216:2) {};

        \draw [edge] (1) -- (2) -- (3) -- (4) (5) -- (1);
        \draw [edge, highlight] (a) -- (1) (b) -- (2) (c) -- (3) (4) -- (5);
        \end{tikzpicture}
        \caption{A graph that contains a $4$-matching.}
        \label{fig:4-matching}
        
    \end{subfigure} \\
    \begin{subfigure}[b]{0.4\textwidth}
        \centering
        
        \begin{tikzpicture}
        \pgftransformrotate{90}
        \node (4) [vertex, label=below left:$4$] at (144:1) {};
        \node (5) [vertex, label=left:$5$] at (72:1) {};
        \node (1) [vertex, label=above:$1$] at (0:1) {};
        \node (2) [vertex, label=right:$2$] at (288:1) {};
        \node (3) [vertex, label=above left:$3$] at (216:1) {};

        \node (a) [vertex, label=left:$a$] at (0:0) {};
        \node (b) [vertex, label=right:$b$] at (216:2) {};

        \draw [edge] (1) -- (2) -- (3) -- (4) -- (5) -- (1);
        \draw [edge] (a) -- (1);
        \draw [edge] (a) -- (2);
        \draw [edge] (b) -- (3);
        \end{tikzpicture}
        \caption{The additional edges off $1$ and $2$ share a vertex that is distinct from the endpoint of the additional edge off $3$.}
        \label{fig:C5CaseI}
        
    \end{subfigure} \qquad
    \begin{subfigure}[b]{0.4\textwidth}
        \centering
        
        \begin{tikzpicture}
        \pgftransformrotate{90}
        \node (4) [vertex, label=below left:$4$] at (144:1) {};
        \node (5) [vertex, label=left:$5$] at (72:1) {};
        \node (1) [vertex, label=above:$1$] at (0:1) {};
        \node (2) [vertex, label=below right:$2$] at (288:1) {};
        \node (3) [vertex, label=below right:$3$] at (216:1) {};

        \node (a) [vertex, label=left:$a$] at (0:0) {};
        \node (b) [vertex, label=above:$b$] at (288:2) {};

        \draw [edge] (1) -- (2) -- (3) -- (4) -- (5) -- (1);
        \draw [edge] (a) -- (1);
        \draw [edge] (b) -- (2);
        \draw [edge] (a) -- (3);
        \end{tikzpicture}
        \caption{The additional edges off $1$ and $3$ share a vertex that is distinct from the endpoint of the additional edge off $2$.}
        \label{fig:C5CaseII}
        
    \end{subfigure}

    \caption{Graphs appearing in the proof of \cref{prop:C_5-implies-C_7}.}
\end{figure}

Now, $N_e$ must be in $\mathcal{G}_3$, so all other edges of $G$ must be either adjacent to $e$ or have both vertices in this $5$-cycle. Since $G$ is connected, we assume the edge $x=1a$ exists without loss of generality. Considering $N_x$, we see that there must be an edge between $b$ and some vertex of $C_5$ other than vertex $1$. If either of the edges $2b$ or $5b$ exist, then $G$ will contain a $7$-cycle. Without loss of generality, we assume that $y=3b$ exists. We see that $N_{45} \in \mathcal{G}_3$, thus $G$ cannot contain any additional vertices.

Observe that adding $4a$, $5b$, or any additional edge with vertex $2$ as an endpoint will create a $7$-cycle. Considering $N_x$ again, we see that it is impossible for this subgraph to contain a path of length four without creating a $7$-cycle in $G$. Therefore, whenever $G$ has an edge disjoint from the $5$-cycle, then $G$ contains a $7$-cycle.

We now instead assume that all other edges of $G$ share at least one vertex with this $5$-cycle. Observe, for example, that the edges $23$ and $45$ form a matching. Since $M(G)$ is two-dimensional and Buchsbaum, these two edges must be part of a $3$-matching with some additional edge whose endpoints are vertex $1$ and some vertex $u_1$ outside of the $5$-cycle. By applying the same argument for each pair of non-adjacent edges in this $5$-cycle, we conclude that each of the five vertices in this $5$-cycle is an endpoint of some edge whose other endpoint is a vertex outside the $5$-cycle. Furthermore, we must be able to choose five such edges---one for each vertex in the $5$-cycle---such that not all share the same new vertex as an endpoint (otherwise any non-adjacent subgraph would be connected but have only four vertices). Thus $G$ has at least two vertices outside this $5$-cycle.

Consider the aforementioned edges for each of the five vertices in the $5$-cycle. If we take three vertices of the $5$-cycle in a row, then at least two of these edges must have a shared endpoint, otherwise $G$ would contain a $4$-matching as depicted in \cref{fig:4-matching}.

Without loss of generality, there are two options for the additional edges off vertices $1$, $2$, and $3$, as depicted in \cref{fig:C5CaseI,fig:C5CaseII}.

Assume that the edges $1a$, $2a$, and $3b$ exist as in \cref{fig:C5CaseI}. Then, considering the vertices $2,3,4$ in a row with the same logic as above, either $4a$ or $4b$ exists. If $4b$ exists, then $G$ contains a $7$-cycle. If $4a$ exists, we apply the same logic to $3,4,5$ to see that either $5a$ or $5b$ exists. If $5b$ exists, then $G$ contains a $7$-cycle. If $5a$ exists, observe that $N_{3b} \in \mathcal{G}_3$, which implies that $G$ cannot have any additional vertices. Furthermore, $N_{23}$ needs a path of length four. The only way to create such a path is to have some edge between one of the vertices in $N_{23}$ and $b$, which creates a $7$-cycle.

Assume instead that the edges $1a$ and $3a$ exist as in \cref{fig:C5CaseII}. We perform the same analysis as in the previous paragraph---either one of the edges $4b$ or $5b$ exists which creates a $7$-cycle, or both of the edges $4a$ and $5a$ exist. In the latter case, we again see that we need some edge between one of the vertices in $N_{23}$ and $b$, which creates a $7$-cycle.
\end{proof}

We know by \cref{prop:C_5-implies-C_7} that the existence of a $5$-cycle will force the existence of a $7$-cycle. We now consider non-bipartite graphs containing $6$-cycles. 

\begin{lem}\label{lem:C_6+C_3}
Let $G$ be a graph whose matching complex is a two-dimensional Buchsbaum complex. If $G$ has both $C_6$ and $C_3$ as subgraphs then it also has $C_7$ as a subgraph.
\end{lem}

\begin{proof}
There are a number of ways that the $C_3$ and $C_6$ could interact. By \cref{prop:C_5-implies-C_7}, any time we deduce that $G$ must have a $C_5$ it must also contain a $C_7$.
\begin{itemize}
\item If the $C_3$ subgraph shares at most one vertex with the $C_6$, then we have a 4-matching---take alternating edges of the $C_6$ together with an edge of the $C_3$ that does not touch the $C_6$. Hence this case cannot appear.
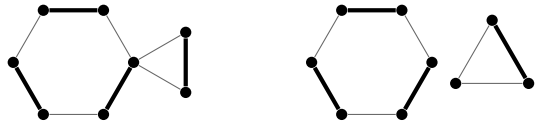
\begin{center}
\begin{tikzpicture}[scale=0.8]
\node (a) [vertex] at (120:1) {};
\node (b) [vertex] at (60:1) {};
\node (c) [vertex] at (0:1) {};
\node (d) [vertex] at (300:1) {};
\node (e) [vertex] at (240:1) {};
\node (f) [vertex] at (180:1) {};

\node (x) [vertex] at ($(c) + (30:1)$) {};
\node (y) [vertex] at ($(c) + (-30:1)$) {};

\draw [edge, gray, thin] (b) -- (c) (d) -- (e) (f) -- (a) (x) -- (c) -- (y);
\draw [edge, ultra thick] (a) -- (b) (c) -- (d) (e) -- (f) (x) -- (y);
\end{tikzpicture} \qquad \qquad
\begin{tikzpicture}[scale=0.8]
\node (a) [vertex] at (120:1) {};
\node (b) [vertex] at (60:1) {};
\node (c) [vertex] at (0:1) {};
\node (d) [vertex] at (300:1) {};
\node (e) [vertex] at (240:1) {};
\node (f) [vertex] at (180:1) {};

\node (x) [vertex] at ($(2,0) + (90:0.7)$) {};
\node (y) [vertex] at ($(2,0) + (210:0.7)$) {};
\node (z) [vertex] at ($(2,0) + (330:0.7)$) {};

\draw [edge, gray, thin] (b) -- (c) (d) -- (e) (f) -- (a) (x) -- (y) -- (z);
\draw [edge, ultra thick] (a) -- (b) (c) -- (d) (e) -- (f) (x) -- (z);
\end{tikzpicture}
\end{center}

\item Suppose our $C_3$ subgraph shares two vertices with the $C_6$ and no edges. These two vertices are either distance $2$ or distance $3$ apart in the $C_6$, and in both cases our graph contains $C_5$ as a subgraph, so \cref{prop:C_5-implies-C_7} gives us a $C_7$.
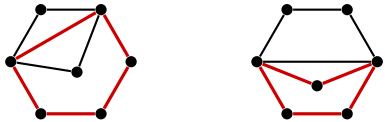
\begin{center}
\begin{tikzpicture}[scale=0.8]
\node (a) [vertex] at (0:1) {};
\node (b) [vertex] at (60:1) {};
\node (c) [vertex] at (120:1) {};
\node (d) [vertex] at (180:1) {};
\node (e) [vertex] at (240:1) {};
\node (f) [vertex] at (300:1) {};

\node (x) [vertex] at (300:0.2) {};

\draw [edge] (b) -- (c) -- (d) -- (x) -- (b);
\draw [edge, highlight] (a) -- (b) -- (d) -- (e) -- (f) -- (a);
\end{tikzpicture} \qquad \qquad
\begin{tikzpicture}[scale=0.8]
\node (a) [vertex] at (0:1) {};
\node (b) [vertex] at (60:1) {};
\node (c) [vertex] at (120:1) {};
\node (d) [vertex] at (180:1) {};
\node (e) [vertex] at (240:1) {};
\node (f) [vertex] at (300:1) {};

\node (x) [vertex] at (270:0.4) {};

\draw [edge] (a) -- (b) -- (c) -- (d) -- (a);
\draw [edge, highlight] (a) -- (x) -- (d) -- (e) -- (f) -- (a);
\end{tikzpicture}
\end{center}

\item Suppose that our $C_3$ subgraph shares three vertices with the $C_6$ and no edges. There is only one way to do this, without loss of generality, and this way gives us a $C_5$ as a subgraph, so again \cref{prop:C_5-implies-C_7} says that $G$ contains $C_7$ as a subgraph.
\begin{center}
\begin{tikzpicture}[scale=0.8]
\node (a) [vertex] at (0:1) {};
\node (b) [vertex] at (60:1) {};
\node (c) [vertex] at (120:1) {};
\node (d) [vertex] at (180:1) {};
\node (e) [vertex] at (240:1) {};
\node (f) [vertex] at (300:1) {};

\draw [edge] (a) -- (b) -- (c) -- (e) -- (a);
\draw [edge, highlight] (a) -- (c) -- (d) -- (e) -- (f) -- (a);
\end{tikzpicture}
\end{center}

\item If our $C_3$ subgraph shares two vertices and one edge with the $C_6$, we immediately get $C_7$ as a subgraph.
\begin{center}
\begin{tikzpicture}[scale=0.8]
\node (a) [vertex] at (0:1) {};
\node (b) [vertex] at (60:1) {};
\node (c) [vertex] at (120:1) {};
\node (d) [vertex] at (180:1) {};
\node (e) [vertex] at (240:1) {};
\node (f) [vertex] at (300:1) {};

\node (x) [vertex] at ($(b) + (1,0)$) {};

\draw [edge] (a) -- (b);
\draw [edge, highlight] (a) -- (x) -- (b) -- (c) -- (d) -- (e) -- (f) -- (a);
\end{tikzpicture}
\end{center}

\item If the $C_3$ shares three vertices and one edge with the $C_6$, we again get $C_5$ as a subgraph, so again we must also have $C_7$.
\begin{center}
\begin{tikzpicture}[scale=0.8]
\node (a) [vertex] at (0:1) {};
\node (b) [vertex] at (60:1) {};
\node (c) [vertex] at (120:1) {};
\node (d) [vertex] at (180:1) {};
\node (e) [vertex] at (240:1) {};
\node (f) [vertex] at (300:1) {};

\draw [edge] (a) -- (d) -- (c) -- (b);
\draw [edge, highlight] (a) -- (b) -- (d) -- (e) -- (f) -- (a);
\end{tikzpicture}
\end{center}

\item Finally, if the $C_3$ shares all three vertices and two edges with the $C_6$, we obtain a $C_5$ and thus also a $C_7$.
\begin{center}
\begin{tikzpicture}[scale=0.8]
\node (a) [vertex] at (0:1) {};
\node (b) [vertex] at (60:1) {};
\node (c) [vertex] at (120:1) {};
\node (d) [vertex] at (180:1) {};
\node (e) [vertex] at (240:1) {};
\node (f) [vertex] at (300:1) {};

\draw [edge] (a) -- (b) -- (c);
\draw [edge, highlight] (a) -- (c) -- (d) -- (e) -- (f) -- (a);
\end{tikzpicture}
\end{center}
\end{itemize}

This covers all possibilities, so $G$ must always contain $C_7$ as a subgraph.
\end{proof}

Once a graph contains a $C_7$, it must have a very constrained structure as described in the lemma below.

\begin{lem}\label{lem:c7-no-extra-vertices}
If $G$ is a graph (with no isolated vertices) that contains $C_7$ as a subgraph and $M(G)$ is two-dimensional, then $G$ has exactly $7$ vertices.
\end{lem}

\begin{figure}
    \centering
    \begin{subfigure}{0.3\textwidth}
    \centering
    \begin{tikzpicture}
    \foreach \i in {1,...,7}
        \node (v\i) [vertex] at (90-\i*360/7:1) {};
        \node (p) [vertex] at (90-2*360/7:2) {};
        \draw [edge, gray, thin] (v1) -- (v2) -- (v3) (v4) -- (v5) (v6) -- (v7);
        \draw [edge, ultra thick] (v2) -- (p) (v3) -- (v4) (v5) -- (v6) (v7) -- (v1);
    \end{tikzpicture}
    \end{subfigure} \qquad
    \begin{subfigure}{0.3\textwidth}
    \centering
    \begin{tikzpicture}
    \foreach \i in {1,...,7}
        \node (v\i) [vertex] at (90-\i*360/7:1) {};
        \node (e1) [vertex] at (2,0.5) {};
        \node (e2) [vertex] at (2,-0.5) {};
        \draw [edge, gray, thin] (v1) -- (v2) -- (v3) (v4) -- (v5) (v6) -- (v7);
        \draw [edge, ultra thick] (e1) -- (e2) (v3) -- (v4) (v5) -- (v6) (v7) -- (v1);
    \end{tikzpicture}
    \end{subfigure}
    \caption{Any graph containing $C_7$ and any additional non-isolated vertices will always contain a $4$-matching.}
    \label{fig:c7-plus-edge-gives-4-matching}
\end{figure}
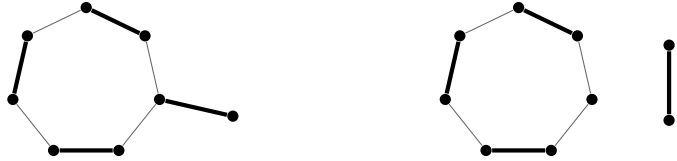

\begin{proof}
See \cref{fig:c7-plus-edge-gives-4-matching}.
\end{proof}

We note that the graphs in \cref{lem:c7-no-extra-vertices} are all Hamiltonian, i.e., they each contain a cycle that uses all vertices of the graph.

A consequence of \cref{lem:c7-no-extra-vertices} is that there are only finitely many possibilities to check to find all graphs containing $C_7$ whose matching complex is two-dimensional and Buchsbaum: Simply take a $C_7$ and add every subset of the $\binom{7}{2} - 7 = 14$ edges that could be added, giving $2^{14} = 16384$ possibilities. While this would be impractical to check by hand, a computer can search these possibilities without difficulty. We have included this code in an appendix. Out of the 383 isomorphism classes of graphs on $7$ vertices containing a $C_7$, 125 of them have a matching complex that is two-dimensional and Buchsbaum. See \cref{tab:c7-data} for more refined data.

\begin{table}
    \centering
    \begin{tabular}{m{3.5cm} | *{15}{c} | c}
        \# edges added to $C_7$ & 0 & 1 & 2 & 3 & 4 & 5 & 6 & 7 & 8 & 9 & 10 & 11 & 12 & 13 & 14 & total \\[1ex]
        \hline
        \hline
        \# graphs up to isomorphism & 1 & 2 & 10 & 30 & 58 & 77 & 73 & 56 & 37 & 20 & 10 & 5 & 2 & 1 & 1 & 383 \\[1ex]
        \hline
        \# graphs where $M(G)$ is 2D Buchsbaum, up to isomorphism & 1 & 1 & 3 & 7 & 11 & 18 & 19 & 20 & 18 & 12 & 7 & 4 & 2 & 1 & 1 & 125
    \end{tabular}
    \caption{Data on graphs containing $C_7$.}
    \label{tab:c7-data}
\end{table}

In particular, $C_7$ itself has a two-dimensional Buchsbaum matching complex, as does $K_7$. Deleting any one or two edges from $K_7$ gives a two-dimensional Buchsbaum matching complex, and so does one of the two ways of adding a single edge to $C_7$ (up to isomorphism).

\Cref{lem:c7-no-extra-vertices} also puts restrictions on what $N_e$ can be for any edge $e$ in these graphs: There must be exactly five vertices in $N_e$, and there are only a few possibilities allowed by \cref{lem:Disconnected-G-Connected-M(G),lem:connected-graph} with only $5$ vertices.

\begin{lem}\label{lem:C7links}
If $G$ is a graph containing $C_7$ whose matching complex is two-dimensional and Buchsbaum, and $e$ is an edge of $G$, then $N_e$ must be one of the graphs in \cref{fig:5-verts-M(G)-connected}.
\end{lem}

\subsection{Link connected graphs}\label{sec:Link-Connected}

We split up our remaining casework using the following definition.

\begin{defn}\label{def:link_conn}
Let $G$ be a connected graph. We call $G$ \emph{link connected} if $N_e$ is a connected graph for every edge $e \in E(G)$.
\end{defn}

\begin{rem}
The above definition is similar to the definition of a 3-vertex-connected graph. In particular a 3-vertex-connected graph (or a ``3-connected graph'') is a graph such that the removal of any two vertices cannot disconnect the graph. The difference here is that we call $G$ link connected if it is a graph such that the removal of any two \emph{adjacent} vertices cannot disconnect the graph. Furthermore, $3$-connected graphs should not have any isolated vertices after removing the two specified vertices, but our definition of $N_e$ omits isolated vertices by construction.
\end{rem}

In this section we will describe all link connected graphs whose matching complexes are two-dimensional Buchsbaum complexes. We start with some tools that will assist us with this.

\begin{lem}\label{lem:cycles}
Suppose $G$ is a link connected graph. If $M(G)$ is a two-dimensional Buchsbaum complex, then $G$ has $C_k$ as a subgraph for some $k \in \set{4,5,6,7}.$
\end{lem}
\begin{proof}

First note that for $k \ge 8$, $G$ cannot have $C_k$ as a subgraph, since such cycles all contain a 4-matching.

Take any edge $uv \in E(G)$ and consider $N_{uv}$. By \cref{lem:path-4} we know that we have the following as a subgraph of $G$:

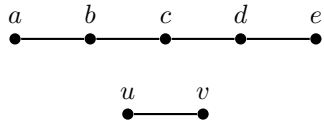
\begin{center}
\begin{tikzpicture}
\node (a) [vertex, label=above:$a$] at (0,0) {};
\node (b) [vertex, label=above:$b$] at (1,0) {};
\node (c) [vertex, label=above:$c$] at (2,0) {};
\node (d) [vertex, label=above:$d$] at (3,0) {};
\node (e) [vertex, label=above:$e$] at (4,0) {};

\node (m) [vertex, label=above:$u$] at (1.5,-1) {};
\node (n) [vertex, label=above:$v$] at (2.5,-1) {};

\draw [edge] (a) -- (b) -- (c) -- (d) -- (e);
\draw [edge] (m) -- (n);
\end{tikzpicture}
\end{center}

Since $G$ is a link connected graph, $N_{bc}$ must be connected, so there must be a path in $G$ between either $u$ or $v$ and either $d$ or $e$. If this path contains $a$, then it creates a $k$-cycle for $k\ge 4$ without using either of $u$ or $v$, so we assume this does not occur. Similarly, $N_{cd}$ must be connected, so there is a path in $G$ from either $u$ or $v$ to either $a$ or $b$. Similarly, we may assume this path does not contain $e$.

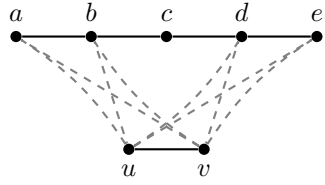
\begin{center}
\begin{tikzpicture}
\node (a) [vertex, label=above:$a$] at (0,0) {};
\node (b) [vertex, label=above:$b$] at (1,0) {};
\node (c) [vertex, label=above:$c$] at (2,0) {};
\node (d) [vertex, label=above:$d$] at (3,0) {};
\node (e) [vertex, label=above:$e$] at (4,0) {};

\node (u) [vertex, label=below:$u$] at (1.5,-1.5) {};
\node (v) [vertex, label=below:$v$] at (2.5,-1.5) {};

\draw [edge] (a) -- (b) -- (c) -- (d) -- (e);
\draw [edge] (u) -- (v);

\draw [edge, optional] (u) to [bend right=12] (a) -- (v) to [bend left=12] (b) -- (u) to [bend right=12] (d) -- (v) to [bend left=12] (e) -- (u);

\end{tikzpicture}
\end{center}

Patching these paths together, we get a cycle whose length is at least $4$.
\end{proof}

We will be using \cref{lem:cycles} to break up our casework in this section.

\begin{lem} \label{lem:C_4-implies-C_567}
Let $G$ be a link connected graph. If $M(G)$ is a two-dimensional Buchsbaum complex and $G$ has $C_4$ as a subgraph, then either $G$ has $C_k$ as a subgraph for some $k \in \set{5,6,7}$ or $G \in \mathcal{B}_1$.
\end{lem}

\begin{proof}
By assumption, $G$ has $C_4$ as a subgraph. Label the vertices of this subgraph as follows:

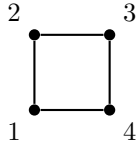
\begin{center}
\begin{tikzpicture}
\node (a) [vertex, label=below left:$1$] at (0,0) {};
\node (b) [vertex, label=above left:$2$] at (0,1) {};
\node (c) [vertex, label=above right:$3$] at (1,1) {};
\node (d) [vertex, label=below right:$4$] at (1,0) {};

\draw [edge] (a) -- (b) -- (c) -- (d) -- (a);
\end{tikzpicture}
\end{center}

Because $12$ and $34$ are in a 2-matching together, they must be part of a 3-matching. Therefore, we must have at least one more edge $56$ which is disjoint from the $4$-cycle:

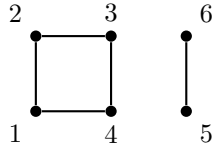
\begin{center}
\begin{tikzpicture}
\node (a) [vertex, label=below left:$1$] at (0,0) {};
\node (b) [vertex, label=above left:$2$] at (0,1) {};
\node (c) [vertex, label=above:$3$] at (1,1) {};
\node (d) [vertex, label=below:$4$] at (1,0) {};
\node (e) [vertex, label=below right:$5$] at (2,0) {};
\node (f) [vertex, label=above right:$6$] at (2,1) {};

\draw [edge] (a) -- (b) -- (c) -- (d) -- (a) (e) -- (f);
\end{tikzpicture}
\end{center}

Since $G$ is link connected, $N_{12}$ must be a connected graph, so without loss of generality there must be a path between vertices $3$ and $6$. This path cannot have length three or more, as that would introduce a $4$-matching, so we must have one of these two cases:

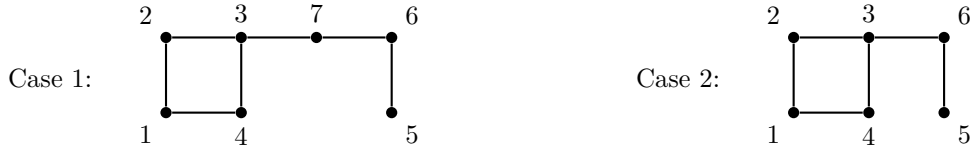
\begin{center}
Case 1: \quad \begin{tikzpicture}
\node (a) [vertex, label=below left:$1$] at (0,0) {};
\node (b) [vertex, label=above left:$2$] at (0,1) {};
\node (c) [vertex, label=above:$3$] at (1,1) {};
\node (d) [vertex, label=below:$4$] at (1,0) {};
\node (e) [vertex, label=below right:$5$] at (3,0) {};
\node (f) [vertex, label=above right:$6$] at (3,1) {};
\node (g) [vertex, label=above:$7$] at (2,1) {};

\draw [edge] (a) -- (b) -- (c) -- (d) -- (a) (c) -- (g) -- (f) -- (e);
\end{tikzpicture} \hspace{1in}
Case 2: \quad \begin{tikzpicture}
\node (a) [vertex, label=below left:$1$] at (0,0) {};
\node (b) [vertex, label=above left:$2$] at (0,1) {};
\node (c) [vertex, label=above:$3$] at (1,1) {};
\node (d) [vertex, label=below:$4$] at (1,0) {};
\node (e) [vertex, label=below right:$5$] at (2,0) {};
\node (f) [vertex, label=above right:$6$] at (2,1) {};

\draw [edge] (a) -- (b) -- (c) -- (d) -- (a) (c) -- (f) -- (e);
\end{tikzpicture}
\end{center}

First, let us consider Case 1. Since $N_{37}$ must be connected, there must be a path from $1$, $2$, or $4$ to either $5$ or $6$, but every way to do this produces a cycle $C_k$ with $k \ge 5$.

This leaves us with Case 2. The edges $14$ and $36$ form a 2-matching together, so they must be part of a 3-matching with some other edge of $G$. This other edge cannot be disjoint from the picture above, since that would give us a $4$-matching made of edges $12$, $34$, $56$ and the new edge. So the possibilities for the new edge are:
\begin{itemize}
\item $25$, in which case $G$ has $C_6$ as a subgraph;
\item $5a$ where $a$ is a vertex not previously in our subgraph, which puts us into Case 1, which we have already dealt with;
\item or $2a$ where $a$ is a vertex not previously in our subgraph.
\end{itemize}

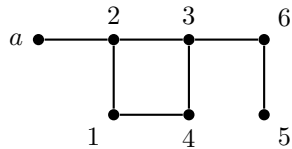
\begin{center}
\begin{tikzpicture}
\node (a) [vertex, label=below left:$1$] at (0,0) {};
\node (b) [vertex, label=above:$2$] at (0,1) {};
\node (c) [vertex, label=above:$3$] at (1,1) {};
\node (d) [vertex, label=below:$4$] at (1,0) {};
\node (e) [vertex, label=below right:$5$] at (2,0) {};
\node (f) [vertex, label=above right:$6$] at (2,1) {};
\node (x) [vertex, label=left:$a$] at (-1,1) {};

\draw [edge] (a) -- (b) -- (c) -- (d) -- (a) (c) -- (f) -- (e) (b) -- (x);
\end{tikzpicture}
\end{center}

In this case, now, $12$ and $36$ are in a 2-matching together, so once again they must be part of a 3-matching. The possibilities for the third edge in this $3$-matching are:
\begin{itemize}
    \item $45$, which gives us a $C_6$;
    \item $4a$, which results in some $C_k$ for $k \in \{5,6,7\}$ after noticing that $N_{23}$ must also be connected;
    \item $4b$ where $b$ is a new vertex, which we will come back to momentarily;
    \item $5a$, which gives a $C_5$;
    \item $5b$, which contains Case 1 and is thus already dealt with;
    \item or $ab$, which gives a $4$-matching.
\end{itemize}
Now we consider the case with edge $4b$.

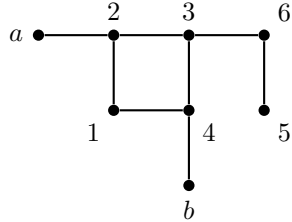
\begin{center}
\begin{tikzpicture}
\node (a) [vertex, label=below left:$1$] at (0,0) {};
\node (b) [vertex, label=above:$2$] at (0,1) {};
\node (c) [vertex, label=above:$3$] at (1,1) {};
\node (d) [vertex, label=below right:$4$] at (1,0) {};
\node (e) [vertex, label=below right:$5$] at (2,0) {};
\node (f) [vertex, label=above right:$6$] at (2,1) {};
\node (x) [vertex, label=left:$a$] at (-1,1) {};
\node (y) [vertex, label=below:$b$] at (1,-1) {};

\draw [edge]
    (a) -- (b) -- (c) -- (d) -- (a)
    (c) -- (f) -- (e)
    (b) -- (x)
    (d) -- (y);
\end{tikzpicture}
\end{center}

Applying similar logic again, $N_{34}$ must be a connected graph, so we can deduce that either $G$ has $C_k$ with $k \in \{5,6,7\}$ or has the following as a subgraph: 

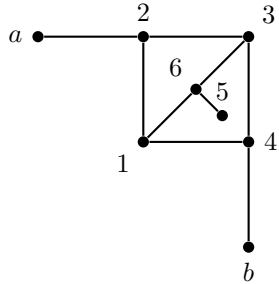
\begin{center}
\begin{tikzpicture}
\pgftransformscale{0.7}
\node (a) [vertex, label=below left:$1$] at (0,0) {};
\node (b) [vertex, label=above:$2$] at (0,2) {};
\node (c) [vertex, label=above right:$3$] at (2,2) {};
\node (d) [vertex, label=right:$4$] at (2,0) {};
\node (e) [vertex, label=above:$5$] at (1.5,0.5) {};
\node (f) [vertex, label=above left:$6$] at (1,1) {};
\node (x) [vertex, label=left:$a$] at (-2,2) {};
\node (y) [vertex, label=below:$b$] at (2,-2) {};

\draw [edge] (a) -- (b) -- (c) -- (d) -- (a)
    (a) -- (f) -- (c)
    (f) -- (e)
    (b) -- (x)
    (d) -- (y);
\end{tikzpicture}
\end{center}

This graph is in $\mathcal{B}_1$ (with vertices $2$, $4$, and $6$ being the unfilled vertices and vertices $1$ and $3$ introduced and attached to all three of the unfilled vertices). If this is a proper subgraph of $G$, then the only other edges we can add without introducing a $4$-matching or a larger cycle are pendants attached to vertices $2$, $4$, and $6$, which keep $G$ in $\mathcal{B}_1$.
\end{proof}

As a consequence, we can split link connected graphs in the following way.

\begin{cor} \label{cor:bipartite-C7-NVC}
Suppose $G$ is link connected and $M(G)$ is a two-dimensional Buchsbaum complex. Then either $G$ has exactly $7$ vertices and has $C_7$ as a subgraph, or $G$ is bipartite.
\end{cor}

\begin{proof}

First, $G$ cannot contain any cycles on $8$ or more vertices without including a $4$-matching.

\Cref{lem:cycles} tells us that $G$ must contain $C_k$, with $4\le k \le 7$. If $G$ includes $C_7$, we are done immediately by \cref{lem:c7-no-extra-vertices}; and if $G$ includes a $C_5$, then \cref{prop:C_5-implies-C_7} gives us a $C_7$ and we are done.

This only leaves the case where $G$ contains $C_4$ or $C_6$ but not $C_5$ or $C_7$. If $G$ contains $C_4$ but not $C_5$ or $C_7$, then \cref{lem:C_4-implies-C_567} implies that either $G$ is in $\mathcal{B}_1$ and thus bipartite, or $G$ contains $C_6$. And if $G$ contains $C_6$ but not $C_5$ or $C_7$, then \cref{lem:C_6+C_3} implies that $G$ cannot contain $C_3$ either, so $G$ has no odd cycles and is thus bipartite. This covers all cases.
\end{proof}

We now consider link connected graphs that contain a $6$-cycle, which will complete our discussion of link connected graphs.

\begin{prop}\label{prop:C_6linkconn}
Suppose $G$ is a link connected graph. If $M(G)$ is a two-dimensional Buchsbaum complex and $G$ has $C_6$ as a subgraph, then $G \in \mathcal{B}_1$, $\mathcal{B}_2$, $\mathcal{B}_3$, $\mathcal{B}_4$, $\mathcal{B}_5$, or $\mathcal{B}_{C_7}$.
\end{prop}

\begin{figure}[b]
    \centering
    \begin{tikzpicture}[scale=1]
        \node (v1) [vertex, label=above:$a$] at (0:0) {};
        \node (g1) [vertex, starpoint, label=right:$1$] at (0:1) {};
        \node (v2) [vertex, label=below right:$2$] at (-60:1) {};
        \node (g2) [vertex, starpoint, label=below left:$3$] at (-120:1) {};
        \node (v3) [vertex, label=left:$4$] at (-180:1) {};
        \node (g3) [vertex, starpoint, label=above left:$5$] at (-240:1) {};
        \node (v4) [vertex, label=above right:$6$] at (-300:1) {};
        \draw [edge] (g1) -- node {$x$} (v2)--(g2) -- node {$y$} (v3)--(g3)--(v4)--(g1) (v1)--(g1);
    \end{tikzpicture}
    \caption{A subgraph appearing in the proof of \cref{prop:C_6linkconn}.}
    \label{fig:C_6linkconn}
\end{figure}
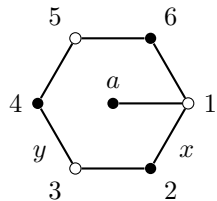

\begin{proof}
Note that if $G$ also contains $C_3$ or $C_5$, then $G \in \mathcal{B}_{C_7}$ by \cref{cor:bipartite-C7-NVC}. Thus we will limit our consideration to bipartite graphs, which implies that $N_e \in \mathcal{G}_1$ for all edges $e$ of $G$. Furthermore, $G$ must contain some edges that have a vertex outside of the $6$-cycle (since each $N_e$ must have at least $5$ vertices), and $G$ cannot have an edge that is disjoint from this $6$-cycle (since this would create a $4$-matching). Therefore we assume $G$ contains the subgraph in \cref{fig:C_6linkconn}. Observe that the filled and unfilled vertices form a bipartition for the graph, so by \cref{prop:Rowan}, there cannot be any more unfilled vertices. (Note that we do not assume that the filled and unfilled vertices follow the convention in \cref{fig:2D-buchsbaum-graphs}; however, these vertices will turn out to follow this convention.)

Considering $N_x$, we see that $G$ must contain an additional edge $e$ with vertex $3$ as an endpoint. Similarly, considering $N_y$, we see that $G$ must contain an additional edge $e'$ with vertex $5$ as an endpoint.

\textbf{Case~1:} The other endpoints of $e$ and $e'$ are not part of the $C_6$.

First we note that if the edges $1a$, $e$, and $e'$ are all pendants, then $G \in \mathcal{B}_1$.

Observe that if $e=3a$ and $e'=5a$, then $G \in \mathcal{B}_3.$ Assume instead $e=3a$ and $e'=5b$ where $b$ is a vertex not shown in \cref{fig:C_6linkconn}. (Observe that this case is equivalent to if $e=3b$ and $e'=5a$ or if $e=3b$ and $e'=5b$.) If neither of the edges $1b$ and $3b$ exist, then $G \in \mathcal{B}_4$.

If both edges $1b$ and $3b$ exist, then again $G \in \mathcal{B}_3$. Assume without loss of generality that the edge $3b$ exists but $1b$ does not. We claim that the family that $G$ is in depends on whether any additional edges between the vertices of the $C_6$ in question exist. If no such edges exist or only the edge $36$ exists, then $G \in \mathcal{B}_2$. If either or both of the other possible edges (i.e., $14$ and $25$) exist, then $G \in \mathcal{B}_3$.

\textbf{Case~2:} At least one of $e$ and $e'$ has both vertices in this $C_6$.

Assume without loss of generality that $e=36$ and consider $N_{36},$ which we recall must be a member of $\mathcal{G}_1.$ Thus $G$ must contain the edge $14$ or one of the edges $25$ and $5a$.

Assume $G$ does not contain either of the edges $25$ or $5a$. Thus $G$ must contain $14$ and $5b$ where $b$ is a vertex not shown in \cref{fig:C_6linkconn}. Therefore $G \in \mathcal{B}_5$.

If instead $G$ contains $5a$, then $G \in \mathcal{B}_3$. Finally, we consider the case where $G$ contains $25$ but not $5a$. Considering $N_{x}$, we see that $G$ must either contain the edge $3a$ (in which case $G \in \mathcal{B}_3$) or either $3b$ or $5b$ where where $b$ is a vertex not shown in \cref{fig:C_6linkconn}. In either of these latter cases, $G \in \mathcal{B}_5$.
\end{proof}

\subsection{Non-link connected graphs}\label{sec:Not-Link-Connected}

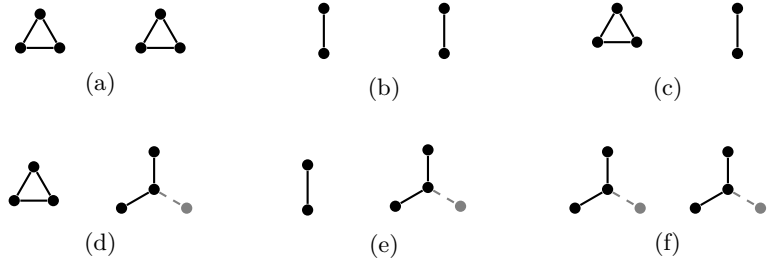
\begin{figure}
    \centering
    \begin{subfigure}{0.2\textwidth}
    \centering
        \begin{tikzpicture}
        \node (a) [coordinate] at (-0.3,0) {};
        \node (b) [coordinate] at (1.3,0) {};
        \node (a1) [vertex] at ($(a)+(90:0.3)$) {};
        \node (a2) [vertex] at ($(a)+(210:0.3)$) {};
        \node (a3) [vertex] at ($(a)+(330:0.3)$) {};
        \node (b1) [vertex] at ($(b)+(90:0.3)$) {};
        \node (b2) [vertex] at ($(b)+(210:0.3)$) {};
        \node (b3) [vertex] at ($(b)+(330:0.3)$) {};
        \draw [edge] (a1) -- (a2) -- (a3) -- (a1) (b1) -- (b2) -- (b3) -- (b1);
        \end{tikzpicture}
        \caption{} \label{fig:nvc-options-a}
    \end{subfigure} \quad 
    \begin{subfigure}{0.2\textwidth}
    \centering
        \begin{tikzpicture}
        \node (a) [coordinate] at (-0.3,0) {};
        \node (b) [coordinate] at (1.3,0) {};
        \node (a1) [vertex] at ($(a)+(0,0.3)$) {};
        \node (a2) [vertex] at ($(a)+(0,-0.3)$) {};
        \node (b1) [vertex] at ($(b)+(0,0.3)$) {};
        \node (b2) [vertex] at ($(b)+(0,-0.3)$) {};
        \draw [edge] (a1) -- (a2) (b1) -- (b2);
        \end{tikzpicture}
        \caption{} \label{fig:nvc-options-b}
    \end{subfigure} \quad 
    \begin{subfigure}{0.2\textwidth}
    \centering
        \begin{tikzpicture}
        \node (a) [coordinate] at (-0.3,0) {};
        \node (b) [coordinate] at (1.3,0) {};
        \node (a1) [vertex] at ($(a)+(90:0.3)$) {};
        \node (a2) [vertex] at ($(a)+(210:0.3)$) {};
        \node (a3) [vertex] at ($(a)+(330:0.3)$) {};
        \node (b1) [vertex] at ($(b)+(0,0.3)$) {};
        \node (b2) [vertex] at ($(b)+(0,-0.3)$) {};
        \draw [edge] (a1) -- (a2) -- (a3) -- (a1) (b1) -- (b2);
        \end{tikzpicture}
        \caption{} \label{fig:nvc-options-c}
    \end{subfigure} \\[0.5cm]
    \begin{subfigure}{0.2\textwidth}
    \centering
        \begin{tikzpicture}
        \node (a) [coordinate] at (-0.3,0) {};
        \node (b) [vertex] at (1.3,0) {};
        \node (a1) [vertex] at ($(a)+(90:0.3)$) {};
        \node (a2) [vertex] at ($(a)+(210:0.3)$) {};
        \node (a3) [vertex] at ($(a)+(330:0.3)$) {};
        \node (b1) [vertex] at ($(b)+(90:0.5)$) {};
        \node (b2) [vertex] at ($(b)+(210:0.5)$) {};
        \node (b3) [vertex, optional] at ($(b)+(330:0.5)$) {};
        \draw [edge] (a1) -- (a2) -- (a3) -- (a1) (b1) -- (b) -- (b2);
        \draw [edge, optional] (b) -- (b3);
        \end{tikzpicture}
        \caption{} \label{fig:nvc-options-e}
    \end{subfigure} \quad
    \begin{subfigure}{0.2\textwidth}
    \centering
        \begin{tikzpicture}
        \node (a) [coordinate] at (-0.3,0) {};
        \node (b) [vertex] at (1.3,0) {};
        \node (a1) [vertex] at ($(a)+(0,0.3)$) {};
        \node (a2) [vertex] at ($(a)+(0,-0.3)$) {};
        \node (b1) [vertex] at ($(b)+(90:0.5)$) {};
        \node (b2) [vertex] at ($(b)+(210:0.5)$) {};
        \node (b3) [vertex, optional] at ($(b)+(330:0.5)$) {};
        \draw [edge] (a1) -- (a2) (b1) -- (b) -- (b2);
        \draw [edge, optional] (b) -- (b3);
        \end{tikzpicture}
        \caption{} \label{fig:nvc-options-d}
    \end{subfigure} \quad
    \begin{subfigure}{0.2\textwidth}
    \centering
        \begin{tikzpicture}
        \node (a) [vertex] at (-0.3,0) {};
        \node (b) [vertex] at (1.3,0) {};
        \node (a1) [vertex] at ($(a)+(90:0.5)$) {};
        \node (a2) [vertex] at ($(a)+(210:0.5)$) {};
        \node (a3) [vertex, optional] at ($(a)+(330:0.5)$) {};
        \node (b1) [vertex] at ($(b)+(90:0.5)$) {};
        \node (b2) [vertex] at ($(b)+(210:0.5)$) {};
        \node (b3) [vertex, optional] at ($(b)+(330:0.5)$) {};
        \draw [edge] (a1) -- (a) -- (a2) (b1) -- (b) -- (b2);
        \draw [edge, optional] (a) -- (a3) (b) -- (b3);
        \end{tikzpicture}
        \caption{} \label{fig:nvc-options-f}
    \end{subfigure}
    \caption{The possibilities for $N_e$ disconnected. \label{fig:nvc-options}}
\end{figure}

We now consider connected graphs which are not link connected. By definition, any such graph must have at least one edge $e$ for which the non-adjacent subgraph $N_e$ is not a connected graph. By \cref{lem:Disconnected-G-Connected-M(G)}, every such connected graph $G$ with a two-dimensional Buchsbaum matching complex must contain an edge $e$ such that $N_e$ is a graph in \cref{fig:nvc-options}. Note for example that under our convention for pendants, each component of the graph in \cref{fig:nvc-options-f} represents any star graph with at least two edges. Though $P_2$ (the path on two vertices) is also a star graph, it does not have a single center vertex and thus behaves somewhat differently from other star graphs, so we consider it separately.

The following lemma shows that any two edges with disconnected non-adjacent subgraphs must share a vertex.

\begin{lem}\label{lem:Multiple-Disconnected}
Let $G$ be a connected graph such that $M(G)$ is two-dimensional. If there exist edges $e,e' \in E$  such that both $N_{e}$ and $N_{e'}$ are disconnected, then $e$ and $e'$ must share a vertex.
\end{lem}

\begin{proof}
Assume that $N_e$ is disconnected with components $G_1$ and $G_2$. All edges of $G$ that are not in one of these components must share a vertex with edge $e$. Without loss of generality, let $e'$ be an edge in $G_1$.

Since $G$ is connected, there must be an edge connecting $e$ and $G_2$. Thus $e$ and $G_2$ form a connected subgraph of $G$, so the only way for $N_{e'}$ to be disconnected is for it to contain some edge $x \in G_1$ that is disjoint from $e'$. But this shows that $\set{e,e',x,y}$ is a $4$-matching for any edge $y \in G_2$, which contradicts that $M(G)$ is two-dimensional.
\end{proof}

We now consider each of the cases from \cref{fig:nvc-options} in turn. As a reminder, $G$ will always be a connected graph in the remainder of this section, and all additional edges of $G$ must be adjacent to the edge $e$.

\begin{prop}\label{prop:nvc(a)}
    Assume $G$ is connected but not link connected and $M(G)$ is a two-dimensional Buchsbaum matching complex. If $G$ contains an edge $e$ such that $N_e = K_3 \sqcup K_3$ (i.e., $N_e$ is the graph in \cref{fig:nvc-options-a}), then $G \in \mathcal{B}_7$.
\end{prop}

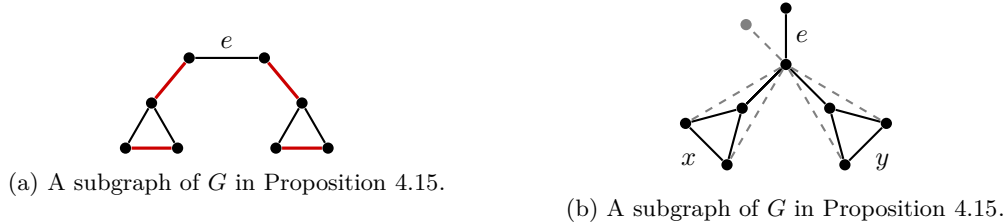
\begin{figure}
    \centering
    \begin{subfigure}{0.4\textwidth}
    \centering
    \begin{tikzpicture}
        \node (e1) [vertex] at (0,1) {};
        \node (e2) [vertex] at (1,1) {};
        \node (a) [coordinate] at (-0.5,0) {};
        \node (b) [coordinate] at (1.5,0) {};
        \node (a1) [vertex] at ($(a)+(90:0.4)$) {};
        \node (a2) [vertex] at ($(a)+(210:0.4)$) {};
        \node (a3) [vertex] at ($(a)+(330:0.4)$) {};
        \node (b1) [vertex] at ($(b)+(90:0.4)$) {};
        \node (b2) [vertex] at ($(b)+(210:0.4)$) {};
        \node (b3) [vertex] at ($(b)+(330:0.4)$) {};
        \draw [edge] (a2) -- (a1) -- (a3) (b2) -- (b1) -- (b3) (e1) -- node {$e$} (e2);
        \draw [edge, highlight] (a2) -- (a3) (a1) -- (e1) (b2) -- (b3) (b1) -- (e2);
        \end{tikzpicture}
    \caption{A subgraph of $G$ in \cref{prop:nvc(a)}.}
    \label{fig:nvc(a)1}
    \end{subfigure} \qquad
    \begin{subfigure}{0.4\textwidth}
    \centering
    \begin{tikzpicture}
    \pgftransformscale{1.5}
    \node (e1) [vertex] at (0.5,1.1) {};
    \node (e2) [vertex] at (0.5,0.6) {};
    \node (es) [vertex, optional] at ($(e2) + (135:0.5)$) {};
    \node (a) [coordinate] at (-0.1,0) {};
    \node (b) [coordinate] at (1.1,0) {};
    \node (a1) [vertex] at ($(a)+(45:0.3)$) {};
    \node (a2) [vertex] at ($(a)+(165:0.3)$) {};
    \node (a3) [vertex] at ($(a)+(285:0.3)$) {};
    \node (b1) [vertex] at ($(b)+(135:0.3)$) {};
    \node (b2) [vertex] at ($(b)+(255:0.3)$) {};
    \node (b3) [vertex] at ($(b)+(15:0.3)$) {};
    \draw [edge] (e1) -- node {$e$} (e2) -- (a1) -- (a2) -- node [swap] {$x$} (a3) -- (a1) (b1) -- (b2) -- node [swap] {$y$} (b3) -- (b1) (e2) -- (a1) (e2) -- (b1);
    \draw [edge, optional] (e2) -- (es) (e2) -- (a2) (e2) -- (a3) (e2) -- (b2) (e2) -- (b3);
    \end{tikzpicture}
    \caption{A subgraph of $G$ in \cref{prop:nvc(a)}.}
    \label{fig:nvc(a)2}
    \end{subfigure}
    \caption{Graphs appearing in \cref{prop:nvc(a)}.}
\end{figure}

\begin{proof}
Since $G$ is connected and all additional edges must touch $e$, each $K_3$ must be connected to e with an edge. If the two copies of $K_3$ are connected to $e$ via different endpoints of $e$, as in \cref{fig:nvc(a)1}, then $G$ will contain a $4$-matching. Since all remaining edges of $G$ must share a vertex with $e$, we therefore must have the graph depicted in \cref{fig:nvc(a)2}, and the other endpoint of $e$ cannot connect to either $K_3$. Considering $N_x$ and $N_y$, we see that each must be a member of $\mathcal{G}_2$. This completely determines $G$, so we see that $G \in \mathcal{B}_7$.
\end{proof}

\begin{prop}\label{prop:nvc(b)}
    Assume $G$ is connected but not link connected and $M(G)$ is a two-dimensional Buchsbaum matching complex. If $G$ contains an edge $e$ such that $N_e = P_2 \sqcup P_2$ (i.e., $N_e$ is the graph in \cref{fig:nvc-options-b}), then $G \in \mathcal{B}_P$.
\end{prop}

\begin{figure}
    \centering
    \begin{subfigure}{0.3\textwidth}
    \centering
    \begin{tikzpicture}
    \pgftransformscale{1.5}
    \node (e1) [vertex] at (0,1) {};
    \node (e2) [vertex] at (1,1) {};
    \node (a) [coordinate] at (-0.5,0) {};
    \node (b) [coordinate] at (1.5,0) {};
    \node (a1) [vertex] at ($(a)+(0,0.3)$) {};
    \node (a2) [vertex] at ($(a)+(0,-0.3)$) {};
    \node (b1) [vertex] at ($(b)+(0,0.3)$) {};
    \node (b2) [vertex] at ($(b)+(0,-0.3)$) {};
    \draw [edge] (a2) -- (a1) to node [auto] {$x$} (e1) to node [auto] {$e$} (e2) to node [auto] {$y$} (b1) -- (b2);
    \end{tikzpicture}
    \caption{A forbidden subgraph of $G$.}
    \label{fig:nvc(b)1}
    \end{subfigure} \qquad
    \begin{subfigure}{0.3\textwidth}
    \centering
    \begin{tikzpicture}
    \pgftransformscale{1.5}
    \node (e1) [vertex, label=right:$a$] at (0.5,0.7) {};
    \node (e2) [vertex, label=right:$b$] at (0.5,1.3) {};
    \node (a) [coordinate] at (-0.1,0) {};
    \node (b) [coordinate] at (1.1,0) {};
    \node (a1) [vertex] at ($(a)+(0,0.3)$) {};
    \node (a2) [vertex] at ($(a)+(0,-0.3)$) {};
    \node (b1) [vertex] at ($(b)+(0,0.3)$) {};
    \node (b2) [vertex] at ($(b)+(0,-0.3)$) {};
    \draw [edge] (a2) to node [auto] {$z$} (a1) -- (e1) to node [auto] {$e$} (e2) (e1) -- (b1) -- (b2);
    \end{tikzpicture}
    \caption{A subgraph of $G$.}
    \label{fig:nvc(b)2}
    \end{subfigure}
    \caption{Graphs appearing in the proof of \cref{prop:nvc(b)}}
\end{figure}
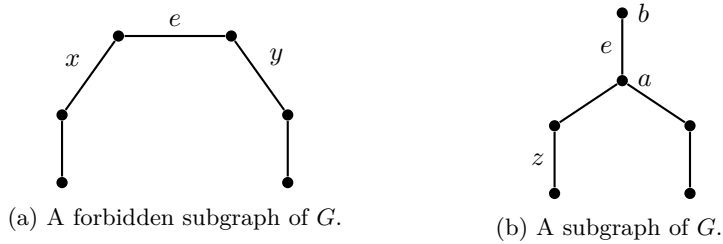

\begin{proof}
Assume that the two copies of $P_2$ are attached to $e$ via different endpoints of $e$. Call these two new edges $x$ and $y$ as in \cref{fig:nvc(b)1}. Since all additional edges of $G$ are adjacent to $e$, there is no 3-matching containing $x$ and $y$, thus $M(G)$ is not Buchsbaum. Thus, only one endpoint of $e$ can be connected to the edges in $N_e$.

This gives the graph depicted in \cref{fig:nvc(b)2}. Since $N_z$ must contain a path of length four, we must add an edge off vertex $b$. Pendants off vertex $a$ would form a $4$-matching and are thus not allowed. We may add any number of pendants off $b$, or, instead, we can add a single edge between $a$ and the pendant edge connected to $b$ to form a $K_3$. Similarly, we may connect up $a$ to either edge in $N_e$. Thus $G \in \mathcal{B}_P$.
\end{proof}

\begin{prop}\label{prop:nvc(c)}
    Assume $G$ is connected but not link connected and $M(G)$ is a two-dimensional Buchsbaum matching complex. If $G$ contains an edge $e$ such that $N_e = K_3 \sqcup P_2$ (i.e., $N_e$ is the graph in \cref{fig:nvc-options-c}), then $G \in \mathcal{B}_7,\mathcal{B}_8,\mathcal{B}_{C_7},$ or $G$ is one of the two exceptional graphs $E_1$ and $E_2$.
\end{prop}

\begin{figure}
    \centering
    \begin{subfigure}[b]{0.4\textwidth}
        \centering
        \begin{tikzpicture}
        \node (bowtie) [coordinate] at (0,0) {};
        \node (edge) [coordinate] at (2.4,0) {};
        \node (e1) [vertex] at ($(bowtie) + (30:1)$) {};
        \node (e2) [vertex] at ($(bowtie) + (-30:1)$) {};
        \node (y1) [vertex] at ($(bowtie) + (0,0)$) {};
        \node (y2) [vertex] at ($(bowtie) + (150:1)$) {};
        \node (v) [vertex] at ($(bowtie) + (210:1)$) {};
        \node (x1) [vertex] at ($(edge) + (0,0.5)$) {};
        \node (x2) [vertex] at ($(edge) + (0,-0.5)$) {};
        \draw [edge] (y1) -- node [auto, swap, xshift=-2pt] {$y$} (y2) -- (v) -- (y1) -- (e1) -- node [auto] {$e$} (e2) -- (y1) (x1) -- node [auto] {$x$} (x2);
        \draw [edge, optional] (e1) to [bend right=20] (x1) -- (e2) to [bend left=20] (x2) -- (e1);
        \end{tikzpicture}
        \caption{The case when $N_x=B$ in \cref{prop:nvc(c)}. We further consider $N_y$, which leads to a contradiction.}
        \label{fig:nvc(c)1}
    \end{subfigure} \\
    \begin{subfigure}[b]{0.4\textwidth}
        \centering
        \begin{tikzpicture}[auto]
        \node (t) [vertex] at (0,0) {};
        \node (y1) [vertex] at (240:1) {};
        \node (y2) [vertex] at (300:1) {};
        \node (e1) [vertex, label=left:$a$] at (0,1) {};
        \node (e2) [vertex, label=left:$b$] at (0,2) {};
        \node (opt) [vertex, optional] at ($(e1) + (135:1)$) {};
        \node (x1) [vertex] at (1.5,2) {};
        \node (x2) [vertex] at (1.5,1) {};
        \draw [edge] (t) -- (y1) -- node [swap] {$y$} (y2) -- (t) -- (e1) -- node [swap] {$e$} (e2) (x1) -- node {$x$} (x2);
        \draw [edge, optional] (e1) -- (y1) (e1) -- (y2) (e1) -- (opt);
        \end{tikzpicture}
        \caption{The case when $N_x \in \mathcal{G}_2$ in \cref{prop:nvc(c)}. Considering $N_y$, we see that there must be an edge connecting $b$ and an endpoint of $x$.}
        \label{fig:nvc(c)2}
    \end{subfigure} \qquad
    \begin{subfigure}[b]{0.4\textwidth}
        \centering
        \begin{tikzpicture}[auto]
        \node (t) [vertex] at (0,0) {};
        \node (y1) [vertex] at (240:1) {};
        \node (y2) [vertex] at (300:1) {};
        \node (e1) [vertex, label=left:$a$] at (0,1) {};
        \node (e2) [vertex, label=above left:$b$] at (0,2) {};
        \node (opt) [vertex, optional] at ($(e1) + (225:1)$) {};
        \node (p) [vertex, highlight] at ($(e1) + (135:1)$) {};
        \node (x1) [vertex] at (1,2) {};
        \node (x2) [vertex] at (1,1) {};
        \draw [edge] (t) -- (y1) -- (y2) -- (t) -- (e1) -- node [swap] {$e$} (e2) -- node {$z$} (x1) -- node {$x$} (x2);
        \draw [edge, optional] (e1) -- (y1) (e1) -- (y2) (e1) -- (opt);
        \draw [edge, highlight] (x2) -- node {$q$} (e1) -- node [pos=0.6, label={[label distance=-4pt]above right:$p$}] {} (p);
        \end{tikzpicture}
        \caption{The case when $N_x \in \mathcal{G}_2$ in \cref{prop:nvc(c)}. Considering $N_z$, least one of the red edges $p$ and $q$ must exist.}
        \label{fig:nvc(c)3}
    \end{subfigure}
    \caption{Graphs appearing in the proof of \cref{prop:nvc(c)}.}
\end{figure}
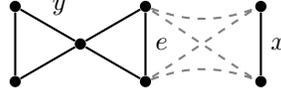
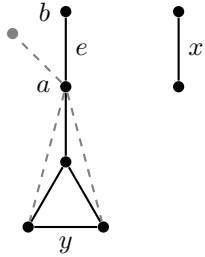
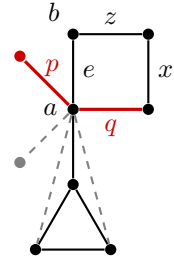

\begin{proof}
Let $x$ be the isolated edge of $N_e$. By \cref{lem:Multiple-Disconnected}, $N_x$ is connected and thus either $N_x=B$ or $N_x \in \mathcal{G}_2~\text{or}~\mathcal{G}_3$. If $N_x \in \mathcal{G}_3$, then $G \in \mathcal{B}_{C_7}$ by \cref{prop:C_5-implies-C_7}.

If $N_x =B$, then the graph in \cref{fig:nvc(c)1} is a subgraph of $G$ and all other edges of $G$ are adjacent to both $e$ and $x$. By \cref{lem:Multiple-Disconnected}, $N_y$ must be connected and thus have a path of length four. But this is impossible, so this case cannot happen.

If instead $N_x \in \mathcal{G}_2$, then the graph in \cref{fig:nvc(c)2} is a subgraph of $G$ and all other edges of $G$ are adjacent to both $e$ and $x$ or are of the form of the dotted edges in \cref{fig:nvc(c)2}. Again, $N_y$ must be connected and have a path of length four, which forces the edge connecting vertex $b$ and an endpoint of $x$ to exist. Calling this edge $z$ and considering $N_z$, at least one of the edges labeled $p$ and $q$ in \cref{fig:nvc(c)3} must appear.

Assume the edge labeled $p$ exists. If the edge between $b$ and the other endpoint of $x$ exists, then $G \in \mathcal{B}_7.$ In not, the only possibility is to have $G \in \mathcal{B}_8.$

If the edge labeled $p$ does not exist, then $q$ must exist. Considering the other edges of the $K_3$ in $N_z$ in turn, we see that the $K_4$ containing vertex $a$ and this $K_3$ must be completed because each edge of this $K_3$ needs a path of length four in its non-adjacent subgraph. Furthermore, the other two edges between $e$ and $x$ can either both exist or neither exist. Thus $G$ is either $E_1$ or $E_2$, the exceptional graphs in \cref{fig:2D-buchsbaum-graphs}.
\end{proof}

\begin{prop}\label{prop:nvc(e)}
    Assume $G$ is connected but not link connected and $M(G)$ is a two-dimensional Buchsbaum matching complex. If $G$ contains an edge $e$ such that $N_e = K_3 \sqcup S_n$ with $n\ge 2$ (i.e., $N_e$ is the graph in \cref{fig:nvc-options-e}), then $G \in \mathcal{B}_8$.
\end{prop}

\begin{figure}
    \centering
    \begin{subfigure}[b]{0.4\textwidth}
        \centering
        \begin{tikzpicture}
        \node (bowtie) [coordinate] at (0,0) {};
        \node (edge) [coordinate] at (2.4,0) {};
        \node (e1) [vertex] at ($(bowtie) + (30:1)$) {};
        \node (e2) [vertex] at ($(bowtie) + (-30:1)$) {};
        \node (y1) [vertex] at ($(bowtie) + (0,0)$) {};
        \node (y2) [vertex] at ($(bowtie) + (150:1)$) {};
        \node (v) [vertex] at ($(bowtie) + (210:1)$) {};
        \node (x1) [vertex] at ($(edge) + (0,0.5)$) {};
        \node (x2) [vertex] at ($(edge) + (0,-0.5)$) {};
        \node (s1) [vertex] at ($(x1) + (120:1)$) {};
        \node (s2) [vertex, optional] at ($(x1) + (60:1)$) {};
        \draw [edge] (y1) -- node [swap] {$y$} (y2) -- (v) -- (y1) -- (e1) -- node {$e$} (e2) -- (y1) (s1) -- (x1) -- node {$x$} (x2);
        \draw [edge, optional] (e1) to [bend right=20] (x1) -- (e2) to [bend left=20] (x2) -- (e1) (x1) -- (s2);
        \end{tikzpicture}
        \caption{The case when $N_x=B$ in \cref{prop:nvc(e)}.}
        \label{fig:nvc(e)1}
    \end{subfigure} \qquad
    \begin{subfigure}[b]{0.4\textwidth}
        \centering
        \begin{tikzpicture}[auto]
        \node (t) [vertex] at (0,0) {};
        \node (y1) [vertex] at (240:1) {};
        \node (y2) [vertex] at (300:1) {};
        \node (e2) [vertex, label=left:$a$] at (0,1) {};
        \node (e1) [vertex, label=left:$b$] at (0,2.5) {};
        \node (opt1) [vertex, optional] at ($(e2) + (150:1)$) {};
        \node (x1) [vertex] at (1.5,2.5) {};
        \node (x2) [vertex] at (1.5,1) {};
        \node (s1) [vertex] at (0.5,2) {};
        \node (s2) [vertex, optional] at ($(x1) + (45:1)$) {};
        \node (opt2) [vertex, optional] at (1,1.5) {};
        \draw [edge] (t) -- (y1) -- node [swap] {$y$} (y2) -- (t) -- (e2) -- node {$e$} (e1) (s1) -- (x1) -- node {$x$} (x2);
        \draw [edge, optional] (e2) -- (y1) (e2) -- (y2) (e2) -- (opt1) (e2) -- (s1) (e2) -- (opt2) -- (x1) -- (s2) (x1) -- (e1) (x2) -- (e2);
        \end{tikzpicture}
        \caption{The case when $N_x \in \mathcal{G}_2$ in \cref{prop:nvc(e)}.}
        \label{fig:nvc(e)2}
    \end{subfigure}
    \caption{Graphs appearing in the proof of \cref{prop:nvc(e)}} \label{fig:nvc(e)}
\end{figure}
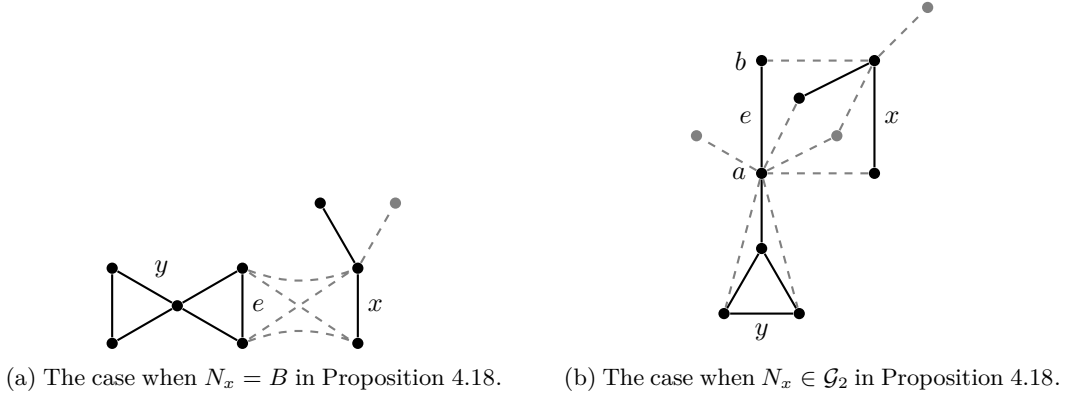

\begin{proof}
Let $x$ be an edge of the star graph in $N_e$ and consider $N_x$. Observe that $N_x \in \mathcal{G}_3$ is impossible by \cref{prop:C_5-implies-C_7,lem:c7-no-extra-vertices}. Since $N_x$ contains a $K_3$, either $N_x = B$ or $N_x \in \mathcal{G}_2$.

Assume $N_x = B$ and we will show that this is impossible. We know that $G$ contains the graph in \cref{fig:nvc(e)1}. Observe that an edge connecting $e$ to a non-central vertex of the star containing $x$ produces a $4$-matching and thus is forbidden. Instead if there is an edge connecting $e$ and the central vertex of the star, we consider $N_y$, where $y$ is the edge labeled in \cref{fig:nvc(e)1}. Observe that $N_y$ is connected but cannot have a path of length four, so this situation is impossible.

Instead assume $N_x \in \mathcal{G}_2$. We claim that $G \in \mathcal{B}_8$. By assumption, $G$ contains the graph \cref{fig:nvc(e)2} as a subgraph, and we claim that any additional edge of $G$ is of the form of a dotted edge in this figure. Observe that we cannot have an edge that connects $b$ to a non-central vertex of the star containing $x$ since this would give a $4$-matching in $G$. Let $y$ be the edge of $K_3$ labeled in \cref{fig:nvc(e)2}. Considering $N_y$, we see that there cannot be an edge that connects $a$ to the central vertex of the star, since this would render $N_y$ connected without a path of length four.

If $G$ has an edge connecting $a$ to a non-central vertex of the star containing $x$, then $G \in \mathcal{B}_8$. If $G$ does not have such an edge, then $G$ must have an edge connecting $b$ to the center of this star. Call this new edge $z$. Considering $N_z$, we see that there must be a pendant off vertex $a$, so again we can conclude that $G \in \mathcal{B}_8$.
\end{proof}

\begin{prop}\label{prop:nvc(d)}
Assume $G$ is connected but not link connected and $M(G)$ is a two-dimensional Buchsbaum matching complex. If $G$ contains an edge $e$ such that $N_e = P_2 \sqcup S_n$ with $n \ge 2$ (i.e., $N_e$ is the graph in \cref{fig:nvc-options-d}), then $G \in \mathcal{B}_2, \mathcal{B}_6, \mathcal{B}_8, \mathcal{B}_9, \mathcal{B}_{C_7}, \text{or}~ \mathcal{B}_P$.
\end{prop}

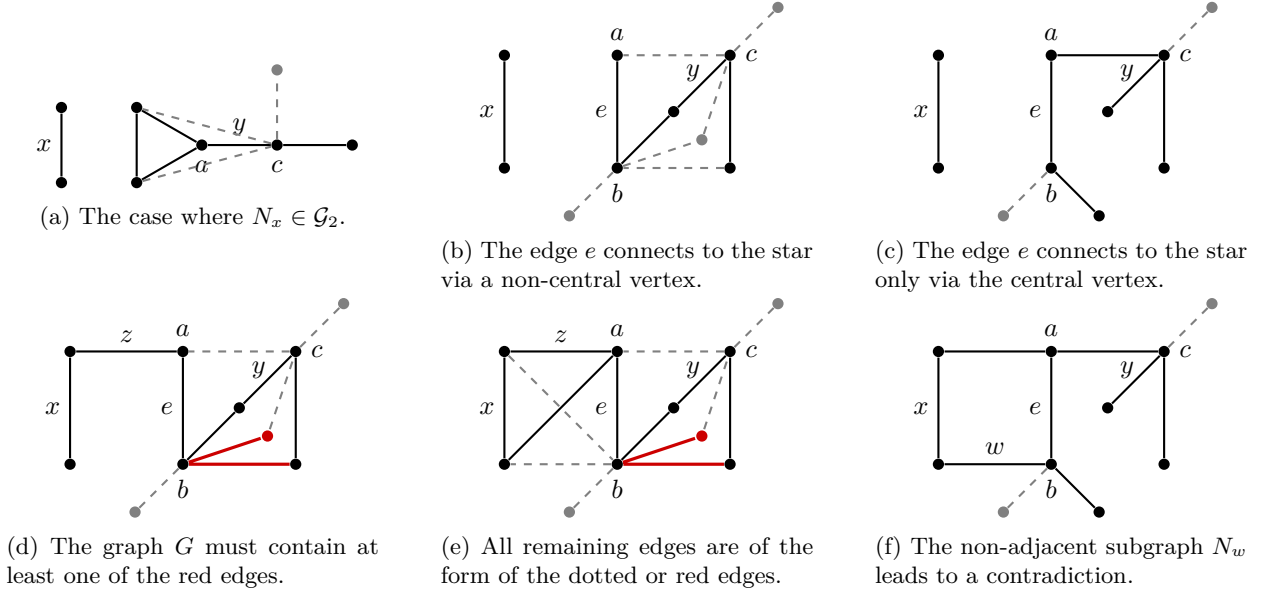
\begin{figure}
    \centering
    \begin{subfigure}{0.3\textwidth}
        \centering
        \begin{tikzpicture}[scale=1]
        \node (x1) [vertex] at (0,1) {};
        \node (x2) [vertex] at (0,0) {};
        \node (e1) [vertex] at (1,1) {};
        \node (e2) [vertex] at (1,0) {};
        \node (y1) [vertex, label=below:$a$] at ($(e2) + (30:1)$) {};
        \node (y2) [vertex, label=below:$c$] at ($(y1) + (1,0)$) {};
        \node (s1) [vertex] at ($(y2) + (1,0)$) {};
        \node (s2) [vertex, optional] at ($(y2) + (0,1)$) {};

        \draw [edge] (x1) -- node [swap] {$x$} (x2) (y1) -- (e1) -- (e2) -- (y1) -- node {$y$} (y2) -- (s1);
        \draw [edge, optional] (e2) -- (y2) -- (e1) (y2) -- (s2);
        \end{tikzpicture}
        \caption{The case where $N_x \in \mathcal{G}_2$.} \label{fig:nvc(d)0}
    \end{subfigure} \qquad
    \begin{subfigure}{0.3\textwidth}
        \centering
        \begin{tikzpicture}[scale=1.5]
        \node (x1) [vertex] at (0,1) {};
        \node (x2) [vertex] at (0,0) {};
        \node (e1) [vertex, label=above:$a$] at (1,1) {};
        \node (e2) [vertex, label=below:$b$] at (1,0) {};
        \node (y1) [vertex] at (1.5,0.5) {};
        \node (y2) [vertex, label=right:$c$] at (2,1) {};
        \node (s1) [vertex] at (2,0) {};
        \node (s2) [vertex, optional] at (1.75,0.25) {};
        \node (s3) [vertex, optional] at ($(y2) + (45:0.6)$) {};
        \node (p) [vertex, optional] at ($(e2) + (225:0.6)$) {};
        
        \draw [edge] (x1) -- node [swap] {$x$} (x2) (e1) -- node [swap] {$e$} (e2) (y1) -- node [shift={(3pt,-3pt)}] {$y$} (y2) -- (s1);
        \draw [edge, optional] (s1) -- (e2) -- (s2) -- (y2) -- (e1) (p) -- (e2) (y2) -- (s3);
        \draw [edge] (e2) -- (y1);
        \end{tikzpicture}
        \caption{The edge $e$ connects to the star via a non-central vertex.} \label{fig:nvc(d)1}
    \end{subfigure} \qquad
    \begin{subfigure}{0.3\textwidth}
        \centering
        \begin{tikzpicture}[scale=1.5]
        \node (x1) [vertex] at (0,1) {};
        \node (x2) [vertex] at (0,0) {};
        \node (e1) [vertex, label=above:$a$] at (1,1) {};
        \node (e2) [vertex, label=below:$b$] at (1,0) {};
        \node (s1) [vertex] at (2,0) {};
        \node (y1) [vertex, label=right:$c$] at (2,1) {};
        \node (y2) [vertex] at (1.5,0.5) {};
        \node (s3) [vertex, optional] at ($(y1) + (45:0.6)$) {};
        \node (p1) [vertex, optional] at ($(e2) + (225:0.6)$) {};
        \node (p2) [vertex] at ($(e2) + (315:0.6)$) {};
        
        \draw [edge] (x1) -- node [swap] {$x$} (x2) (e1) -- node [swap] {$e$} (e2) (s1) -- (y1) -- node [swap, shift={(3pt,-3pt)}] {$y$} (y2);
        \draw [edge, optional] (p1) -- (e2) (y1) -- (s3);
        \draw [edge] (e2) -- (p2) (e1) -- (y1);
        \end{tikzpicture}
        \caption{The edge $e$ connects to the star only via the central vertex.} \label{fig:nvc(d)2}
    \end{subfigure} \qquad
    \begin{subfigure}{0.3\textwidth}
        \centering
        \begin{tikzpicture}[scale=1.5]
        \node (x1) [vertex] at (0,1) {};
        \node (x2) [vertex] at (0,0) {};
        \node (e1) [vertex, label=above:$a$] at (1,1) {};
        \node (e2) [vertex, label=below:$b$] at (1,0) {};
        \node (y1) [vertex] at (1.5,0.5) {};
        \node (y2) [vertex, label=right:$c$] at (2,1) {};
        \node (s1) [vertex] at (2,0) {};
        \node (s2) [vertex, highlight] at (1.75,0.25) {};
        \node (s3) [vertex, optional] at ($(y2) + (45:0.6)$) {};
        \node (p) [vertex, optional] at ($(e2) + (225:0.6)$) {};
        
        \draw [edge] (x2) -- node {$x$} (x1) -- node {$z$} (e1) -- node [swap] {$e$} (e2) -- (y1) -- node [shift={(3pt,-3pt)}] {$y$} (y2) -- (s1);
        \draw [edge, optional] (s2) -- (y2) -- (e1) (p) -- (e2) (y2) -- (s3);
        \draw [edge, highlight] (s1) -- (e2) -- (s2);
        \end{tikzpicture}
        \caption{The graph $G$ must contain at least one of the red edges.} \label{fig:nvc(d)3}
    \end{subfigure} \qquad
    \begin{subfigure}{0.3\textwidth}
        \centering
        \begin{tikzpicture}[scale=1.5]
        \node (x1) [vertex] at (0,1) {};
        \node (x2) [vertex] at (0,0) {};
        \node (e1) [vertex, label=above:$a$] at (1,1) {};
        \node (e2) [vertex, label=below:$b$] at (1,0) {};
        \node (y1) [vertex] at (1.5,0.5) {};
        \node (y2) [vertex, label=right:$c$] at (2,1) {};
        \node (s1) [vertex] at (2,0) {};
        \node (s2) [vertex, highlight] at (1.75,0.25) {};
        \node (s3) [vertex, optional] at ($(y2) + (45:0.6)$) {};
        \node (p) [vertex, optional] at ($(e2) + (225:0.6)$) {};
        
        \draw [edge] (e1) -- (x2) -- node {$x$} (x1) -- node {$z$} (e1) -- node [swap] {$e$} (e2) -- (y1) -- node [shift={(3pt,-3pt)}] {$y$} (y2) -- (s1);
        \draw [edge, optional] (s2) -- (y2) -- (e1) (p) -- (e2) (y2) -- (s3) (x1) -- (e2) -- (x2);
        \draw [edge, highlight] (s1) -- (e2) -- (s2);
        \end{tikzpicture}
        \caption{All remaining edges are of the form of the dotted or red edges.} \label{fig:nvc(d)3.5}
    \end{subfigure} \qquad
    \begin{subfigure}{0.3\textwidth}
        \centering
        \begin{tikzpicture}[scale=1.5]
        \node (x1) [vertex] at (0,1) {};
        \node (x2) [vertex] at (0,0) {};
        \node (e1) [vertex, label=above:$a$] at (1,1) {};
        \node (e2) [vertex, label=below:$b$] at (1,0) {};
        \node (s1) [vertex] at (2,0) {};
        \node (y1) [vertex, label=right:$c$] at (2,1) {};
        \node (y2) [vertex] at (1.5,0.5) {};
        \node (s3) [vertex, optional] at ($(y1) + (45:0.6)$) {};
        \node (p1) [vertex, optional] at ($(e2) + (225:0.6)$) {};
        \node (p2) [vertex] at ($(e2) + (315:0.6)$) {};
        
        \draw [edge] (x1) -- node [swap] {$x$} (x2) -- node {$w$} (e2) -- node {$e$} (e1) -- (x1) (s1) -- (y1) -- node [swap, shift={(3pt,-3pt)}] {$y$} (y2) (e2) -- (p2) (e1) -- (y1);
        \draw [edge, optional] (p1) -- (e2) (s3) -- (y1);
        \end{tikzpicture}
        \caption{The non-adjacent subgraph $N_w$ leads to a contradiction.} \label{fig:nvc(d)4}
    \end{subfigure} \qquad
    \caption{Graphs appearing in the proof of  \cref{prop:nvc(d)}.}
    \label{fig:nvc(d)}
\end{figure}

\begin{proof}
Let $x$ be the isolated edge in $N_e$, let $c$ be the center of the $S_n$ in $N_e$, and consider $N_x$. Since all edges of $N_x$ besides $e$ and the edges in the $S_n$ are adjacent to $e$, $N_x = B$ is impossible, and if $N_x \in \mathcal{G}_3$, then $G \in \mathcal{B}_{C_7}$ by \cref{prop:C_5-implies-C_7}.

First assume that $N_x \in \mathcal{G}_2$. In this case, $G$ contains the graph in \cref{fig:nvc(d)0} as a subgraph, and $e$ is one of the edges of the solid $K_3$ in this figure. All additional edges of $G$ are of the form of one of the dotted edges or are between $e$ and $x$.

If $e$ is the vertical edge in the $K_3$ in \cref{fig:nvc(d)0}, we let $y$ be the edge labeled in this figure and consider $N_y$. By \cref{lem:Multiple-Disconnected}, $N_y$ is connected and thus must have a path of length four. However, this is impossible given the restrictions on $G$; thus we conclude that $e$ is one of the other edges of the solid triangle in this figure. We must add either or both edges between vertex $a$ and edge $x$ to ensure that $G$ is connected, and we may also add pendants off the center of the star $S_n$; no other edges are allowed. In either case, $G \in \mathcal{B}_P$ with one or two petals being $K_3$.

The last case to consider is when $N_x \in \mathcal{G}_1$. We will consider two cases: In the first $e$ connects to the star in $N_e$ via a non-central vertex (\cref{fig:nvc(d)1}); in the second no such connections are allowed (\cref{fig:nvc(d)2}). 
In each case, all remaining edges in $N_x$ are of the form of a dotted edge in the respective figure. All other edges of $G$ must be adjacent to both $e$ and $x$.

\textbf{Case~1:} We have the case in \cref{fig:nvc(d)1}; i.e., the star containing $c$ connects to edge $e$ via at least one non-central vertex.

Consider $N_y$ where $y$ is the edge indicated in this figure. Note that $N_y$ must be connected and thus have a path of length four. Thus $G$ must have an edge connecting edge $x$ and the vertex $a$ (call this edge $z$) and also at least one of the red edges in \cref{fig:nvc(d)3}.

If $N_y \in \mathcal{G}_1$, then the only possible addition to $G$---apart from the edges in \cref{fig:nvc(d)3}---is the edge which connects vertex $b$ to the endpoint of $x$ that is not an endpoint of $z$. 
Assume this edge exists and call it $w$. In this case $G$ cannot contain the edge $ac$, since $N_{ac}$ would be connected without a path of length four. If the red edge connecting $b$ and a non-central vertex of the star exists then $G \in \mathcal{B}_2$, and if this edge does not exist then $G \in \mathcal{B}_6$. If instead $G$ does not contain the edge $w$, then $G \in \mathcal{B}_6$. (To see this, consider $N_z$ in \cref{fig:nvc(d)3} to see we must have an additional edge off $b$ or $c$.)

If instead $N_y \in \mathcal{G}_2$, then $G$ must contain the edge that creates a triangle with edges $x$ and $z$. Thus $G$ contains the graph in \cref{fig:nvc(d)3.5} and all other edges of $G$ are of the form of the dotted or red edges in this figure. Assume the edge between $a$ and $c$ exists. Considering the non-adjacent subgraph for this edge, we see that no additional edges between $e$ and $x$ can exist. Considering $N_z$, we see that either $b$ or $c$ must have a pendant and thus $G \in \mathcal{B}_9$.

Assume instead that the edge between $a$ and $c$ does not exist. If the other two edges between $e$ and $x$ both exist, then either $G$ is the exceptional graph $E_2$ or $G \in \mathcal{B}_8$. If $G$ contains only one of the other possible edges between $e$ and $x$ we assume, without loss of generality, that it is the edge adjacent to $z$. Considering $N_z$ shows $G$ must have another edge off $b$ or $c$ (that doesn't attach to $a$). Thus $G \in \mathcal{B}_8$.

\textbf{Case~2:} We have the case in \cref{fig:nvc(d)2}; i.e., the star containing $c$ is to connected to edge $e$ only via the edge $ac$.

Consider $N_y$ where $y$ is the edge labeled in this figure. For the same reasons as the previous case, there must be an edge between $x$ and vertex $a$, and all other edges of $G$ are either between $e$ and $x$ or of the form of a dotted edge in this figure.

Assume $N_y \in \mathcal{G}_1$. The only way for $N_y$ to contain a cycle in this case is to have edge $w$ as in \cref{fig:nvc(d)4}. However, we see that $N_w$ cannot contain a path of length four, so this edge cannot exist. Thus in this case $G$ is the graph in \cref{fig:nvc(d)2} with an additional edge from $x$ to vertex $a$. So $G \in \mathcal{B}_P$ with all three petals being star graphs.

If instead $N_y \in \mathcal{G}_2$, then the edge between vertex $a$ and the other endpoint of $x$ must exist. Furthermore, neither of the other two edges connecting $e$ and $x$ can exist. Thus $G \in \mathcal{B}_P$ with one petal being a $K_3$. This completes the proof.
\end{proof}

\begin{prop}\label{prop:nvc(f)}
    Assume $G$ is connected but not link connected and $M(G)$ is a two-dimensional Buchsbaum matching complex. If $G$ contains an edge $e$ such that $N_e = S_m \sqcup S_n$ with $m,n\ge 2$ (i.e., $N_e$ is the graph in \cref{fig:nvc-options-f}), then $G \in \mathcal{B}_2, \mathcal{B}_5, \mathcal{B}_6, \text{or}~\mathcal{B}_P$.
\end{prop}

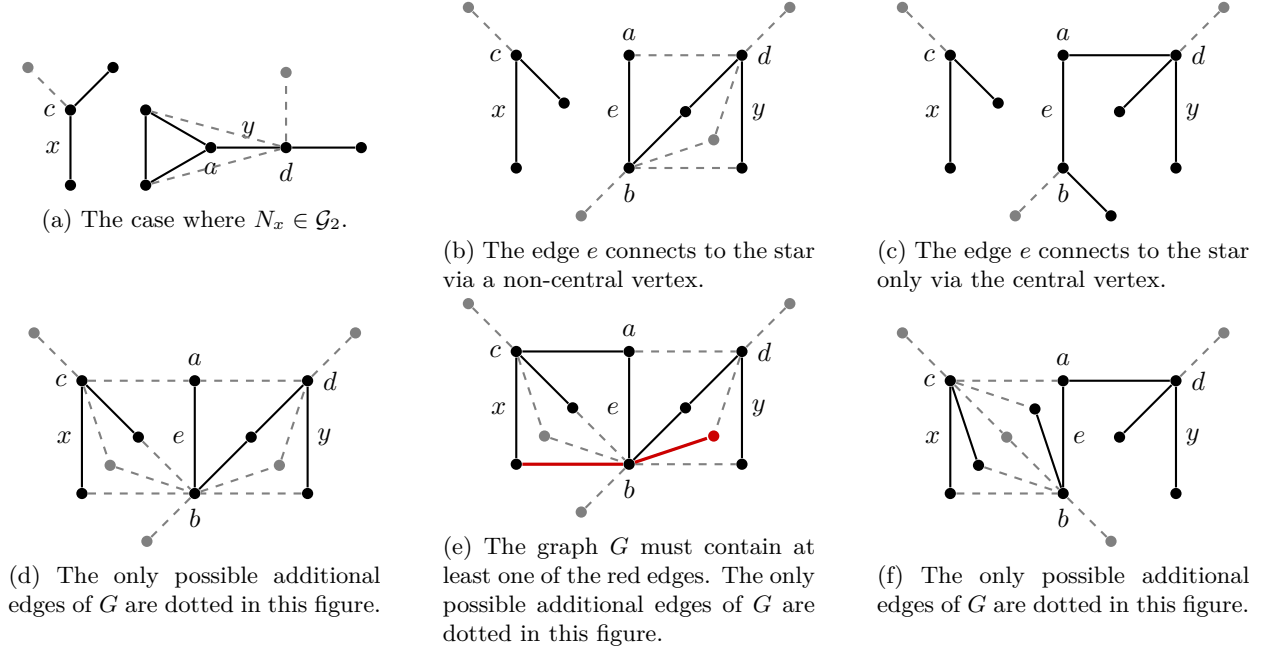
\begin{figure}
    \centering
    \begin{subfigure}{0.3\textwidth}
        \centering
        \begin{tikzpicture}[scale=1]
        \node (x1) [vertex, label=left:$c$] at (0,1) {};
        \node (x2) [vertex] at (0,0) {};
        \node (e1) [vertex] at (1,1) {};
        \node (e2) [vertex] at (1,0) {};
        \node (y1) [vertex, label=below:$a$] at ($(e2) + (30:1)$) {};
        \node (y2) [vertex, label=below:$d$] at ($(y1) + (1,0)$) {};
        \node (s1) [vertex] at ($(y2) + (1,0)$) {};
        \node (s2) [vertex, optional] at ($(y2) + (0,1)$) {};
        
        \node (z1) [vertex] at ($(x1) + (45:0.8)$) {};
        \node (z2) [vertex, optional] at ($(x1) + (135:0.8)$) {};
        \draw [edge] (x1) -- (z1);
        \draw [edge, optional] (x1) -- (z2);

        \draw [edge] (x1) -- node [swap] {$x$} (x2) (y1) -- (e1) -- (e2) -- (y1) -- node {$y$} (y2) -- (s1);
        \draw [edge, optional] (e2) -- (y2) -- (e1) (y2) -- (s2);
        \end{tikzpicture}
        \caption{The case where $N_x \in \mathcal{G}_2$.} \label{fig:nvc(f)1}
    \end{subfigure} \qquad
    \begin{subfigure}{0.3\textwidth}
        \centering
        \begin{tikzpicture}[scale=1.5]
        \node (x1) [vertex, label=left:$c$] at (0,1) {};
        \node (x2) [vertex] at (0,0) {};
        \node (e1) [vertex, label=above:$a$] at (1,1) {};
        \node (e2) [vertex, label=below:$b$] at (1,0) {};
        \node (y1) [vertex] at (1.5,0.5) {};
        \node (y2) [vertex, label=right:$d$] at (2,1) {};
        \node (s1) [vertex] at (2,0) {};
        \node (s2) [vertex, optional] at (1.75,0.25) {};
        \node (s3) [vertex, optional] at ($(y2) + (45:0.6)$) {};
        \node (p) [vertex, optional] at ($(e2) + (225:0.6)$) {};
        
        \node (z1) [vertex] at ($(x1) + (-45:0.6)$) {};
        \node (z2) [vertex, optional] at ($(x1) + (135:0.6)$) {};
        \draw [edge] (x1) -- (z1);
        \draw [edge, optional] (x1) -- (z2);
        
        \draw [edge] (x1) -- node [swap] {$x$} (x2) (e1) -- node [swap] {$e$} (e2) (y1) -- (y2) -- node  {$y$}  (s1);
        \draw [edge, optional] (s1) -- (e2) -- (s2) -- (y2) -- (e1) (p) -- (e2) (y2) -- (s3);
        \draw [edge] (e2) -- (y1);
        \end{tikzpicture}
        \caption{The edge $e$ connects to the star via a non-central vertex.} \label{fig:nvc(f)2}
    \end{subfigure} \qquad
    \begin{subfigure}{0.3\textwidth}
        \centering
        \begin{tikzpicture}[scale=1.5]
        \node (x1) [vertex, label=left:$c$] at (0,1) {};
        \node (x2) [vertex] at (0,0) {};
        \node (e1) [vertex, label=above:$a$] at (1,1) {};
        \node (e2) [vertex, label=below:$b$] at (1,0) {};
        \node (s1) [vertex] at (1.5,0.5) {};
        \node (y1) [vertex, label=right:$d$] at (2,1) {};
        \node (y2) [vertex] at (2,0) {};
        \node (s3) [vertex, optional] at ($(y1) + (45:0.6)$) {};
        \node (p1) [vertex, optional] at ($(e2) + (225:0.6)$) {};
        \node (p2) [vertex] at ($(e2) + (315:0.6)$) {};

        \node (z1) [vertex] at ($(x1) + (-45:0.6)$) {};
        \node (z2) [vertex, optional] at ($(x1) + (135:0.6)$) {};
        \draw [edge] (x1) -- (z1);
        \draw [edge, optional] (x1) -- (z2);
        
        \draw [edge] (x1) -- node [swap] {$x$} (x2) (e1) -- node [swap] {$e$} (e2) (s1) -- (y1) -- node {$y$} (y2);
        \draw [edge, optional] (p1) -- (e2) (y1) -- (s3);
        \draw [edge] (e2) -- (p2) (e1) -- (y1);
        \end{tikzpicture}
        \caption{The edge $e$ connects to the star only via the central vertex.} \label{fig:nvc(f)3}
    \end{subfigure} \qquad
    \begin{subfigure}{0.3\textwidth}
        \centering
        \begin{tikzpicture}[scale=1.5]
        \node (x1) [vertex, label=left:$c$] at (0,1) {};
        \node (x2) [vertex] at (0,0) {};
        \node (e1) [vertex, label=above:$a$] at (1,1) {};
        \node (e2) [vertex, label=below:$b$] at (1,0) {};
        \node (y1) [vertex] at (1.5,0.5) {};
        \node (y2) [vertex, label=right:$d$] at (2,1) {};
        \node (s1) [vertex] at (2,0) {};
        \node (s2) [vertex, optional] at (1.75,0.25) {};
        \node (s3) [vertex, optional] at ($(y2) + (45:0.6)$) {};
        \node (p) [vertex, optional] at ($(e2) + (225:0.6)$) {};

        \node (z1) [vertex] at (0.5,0.5) {};
        \node (z2) [vertex, optional] at ($(x1) + (135:0.6)$) {};
        \node (z3) [vertex, optional] at (0.25,0.25) {};
        \draw [edge] (x1) -- (z1);
        \draw [edge, optional] (x1) -- (z2);
        \draw [edge, optional] (z1) -- (e2) (x1) -- (e1) (x2) -- (e2) (x1) -- (z3) -- (e2);

        \draw [edge] (x2) -- node {$x$} (x1) (e1) -- node [swap] {$e$} (e2) -- (y1) -- (y2) -- node {$y$} (s1);
        \draw [edge, optional] (s2) -- (y2) -- (e1) (p) -- (e2) (y2) -- (s3) (s1) -- (e2) -- (s2);
        \end{tikzpicture}
        \caption{The only possible additional edges of $G$ are dotted in this figure.} \label{fig:nvc(f)4}
    \end{subfigure} \qquad
    \begin{subfigure}{0.3\textwidth}
        \centering
        \begin{tikzpicture}[scale=1.5]
        \node (x1) [vertex, label=left:$c$] at (0,1) {};
        \node (x2) [vertex] at (0,0) {};
        \node (e1) [vertex, label=above:$a$] at (1,1) {};
        \node (e2) [vertex, label=below:$b$] at (1,0) {};
        \node (s1) [vertex] at (1.5,0.5) {};
        \node (s2) [vertex, highlight] at (1.75,0.25) {};
        \node (y1) [vertex, label=right:$d$] at (2,1) {};
        \node (y2) [vertex] at (2,0) {};
        \node (s3) [vertex, optional] at ($(y1) + (45:0.6)$) {};
        \node (p1) [vertex, optional] at ($(e2) + (225:0.6)$) {};

        \node (z1) [vertex] at (0.5,0.5) {};
        \node (z2) [vertex, optional] at ($(x1) + (135:0.6)$) {};
        \node (z3) [vertex, optional] at (0.25,0.25) {};
        \draw [edge] (x1) -- (z1);
        \draw [edge, optional] (x1) -- (z2);
        \draw [edge, optional] (z1) -- (e2) (x1) -- (z3) -- (e2); 
     
        \draw [edge] (x1) -- node [swap] {$x$} (x2) (e2) -- node {$e$} (e1) -- (x1) (e2) -- (s1) -- (y1) -- node {$y$} (y2);
        \draw [edge, optional] (p1) -- (e2) (s3) -- (y1) (e1) -- (y1) (e2) -- (y2) (y1) -- (s2) -- (e2);
        
        \draw[edge, highlight] (e2) -- (x2) (e2) -- (s2);
        \end{tikzpicture}
        \caption{The graph $G$ must contain at least one of the red edges. The only possible additional edges of $G$ are dotted in this figure.} \label{fig:nvc(f)5}
    \end{subfigure} \qquad
    \begin{subfigure}{0.3\textwidth}
        \centering
        \begin{tikzpicture}[scale=1.5]
        \node (x1) [vertex, label=left:$c$] at (0,1) {};
        \node (x2) [vertex] at (0,0) {};
        \node (e1) [vertex, label=above:$a$] at (1,1) {};
        \node (e2) [vertex, label=below:$b$] at (1,0) {};
        \node (y1) [vertex] at (1.5,0.5) {};
        \node (y2) [vertex, label=right:$d$] at (2,1) {};
        \node (s1) [vertex] at (2,0) {};
        \node (s3) [vertex, optional] at ($(y2) + (45:0.6)$) {};
        \node (p) [vertex, optional] at ($(e2) + (-45:0.6)$) {};

        \node (z1) [vertex, optional] at (0.5,0.5) {};
        \node (z2) [vertex, optional] at ($(x1) + (135:0.6)$) {};
        \node (z3) [vertex] at (0.25,0.25) {};
        \node (z4) [vertex] at (0.75,0.75) {};
        \draw [edge] (x1) -- (z3) (e2) -- (z4);
        \draw [edge, optional] (x1) -- (z1) (x1) -- (z2);
        \draw [edge, optional] (z1) -- (e2) (x1) -- (e1) (x2) -- (e2) (z3) -- (e2) (x1) -- (z4);

        \draw [edge] (x2) -- node {$x$} (x1) (e1) -- node {$e$} (e2) (y1) -- (y2) -- node {$y$} (s1) (e1) -- (y2);
        \draw [edge, optional] (y2) -- (s3) (e2) -- (p);
        \end{tikzpicture}
        \caption{The only possible additional edges of $G$ are dotted in this figure.} \label{fig:nvc(f)6}
    \end{subfigure} \qquad
    \caption{Graphs appearing in the proof of  \cref{prop:nvc(f)}.}
    \label{fig:nvc(f)}
\end{figure}

\begin{proof}
Let $c$ be the center vertex of $S_m$ and $d$ be the center vertex of $S_n$. Let $x$ be any edge in $S_m$ and consider $N_x$. Since all remaining edges of $N_x$ are adjacent to $e$, $N_x = B$ is impossible. Similarly $N_x \in \mathcal{G}_3$ is impossible by \cref{prop:C_5-implies-C_7,lem:c7-no-extra-vertices} since $G$ has at least eight vertices.

Assume that  $N_x \in \mathcal{G}_2$. In this case, $G$ contains the graph in \cref{fig:nvc(f)1} as a subgraph, and $e$ is one of the edges of the solid $K_3$ in this figure. All remaining edges of $G$ are either of the form of the dotted edges in this figure or adjacent to both $e$ and $x$. As in \cref{prop:nvc(d)}, we see that $e$ cannot be the vertical edge in \cref{fig:nvc(f)1} by considering $N_y$. Similarly, if $e$ is a different edge of the triangle, then $G$ must have an edge between vertices $a$ and $c$. Furthermore, the only edges we can add to $G$ are more pendants off these stars (and not the dotted edges completing the $K_4$ in \cref{fig:nvc(f)1}), so $G \in \mathcal{B}_P$ with two pendants being stars and the other a $K_3$.

The only remaining possibility is to have $N_x \in \mathcal{G}_1$. Observe that $e$ can connect to the star with center vertex $d$ via a non-central vertex as in \cref{fig:nvc(f)2} (and possibly also $d$) or only via the central vertex $d$ as in \cref{fig:nvc(f)3}. All remaining edges in $N_x$ have to be of the form of one of the dotted edges in these figures. All other edges of $G$ must be adjacent to both $e$ and $x$ or to the vertex $c$.

\textbf{Case~1:} We have the case in \cref{fig:nvc(f)2}; i.e., the star containing $d$ connects to edge $e$ via at least one non-central vertex.

Consider $N_y$ where $y$ is the edge indicated in this figure. Note that $N_y$ must be connected and thus have a path of length four. 
Observe that connecting $b$ to $c$ or connecting $a$ to a non-central vertex of the star containing $c$ both create a contradiction for $N_y$.
Thus the only remaining edges in $G$ must be of the form of the dotted edges in \cref{fig:nvc(f)4}.

There must be some edge connecting $e$ and the star containing $x$. Assume the edge $ac$ does not exist. Without loss of generality, this implies that the edge connecting $b$ to $x$ exists. If, for each star, there is at least one more edge connecting $e$ to the star, then $G \in \mathcal{B}_2$. If not, then $G \in \mathcal{B}_6$.

Assume instead that the edge $ac$ does exist. Considering $N_{ac}$, this forces the existence of an additional edge containing $b$ that is not adjacent to edge $y$. These are indicated in red in \cref{fig:nvc(f)5}. If there are no edges besides $ac$ connecting $e$ and the star containing $c$, then $G \in \mathcal{B}_6$.  If instead there exists an additional edge connecting $b$ to the star containing $c$, we consider a few options. If the edge $ad$ exists, then $G \in \mathcal{B}_5$. If $ad$ does not exist but the edge between $b$ and $y$ exists, then $G \in \mathcal{B}_2$. If neither $ad$ nor the edge between $b$ and $y$ exist, then $G \in \mathcal{B}_6$.

\textbf{Case~2:} We have the case in \cref{fig:nvc(f)3}; i.e., the star containing $d$ is to connected to edge $e$ only via the edge $ad$.

Let $y$ be the edge indicated in this \cref{fig:nvc(f)3}. 
Considering $N_y$, we see that adding either the edge from $b$ to $c$ or the edge from $a$ to a non-central vertex of the star containing $c$ would create a contradiction.

Thus all additional edges in $G$ must be of the form of the dotted edges in \cref{fig:nvc(f)6}. There must be some edge connecting $e$ and the star containing $x$. If the edge between $c$ and the end of the pendant off vertex $b$ in \cref{fig:nvc(f)6} exists, then an additional edge off $b$ must exist, which shows that $G \in \mathcal{B}_6$. Similarly, if a different edge exists between $b$ and a non-central vertex of this star, then $G \in \mathcal{B}_6$. If instead the only edge connecting $e$ and the star containing $x$ is the edge between vertices $a$ and $c$, then $G \in \mathcal{B}_P$ (with all three petals being stars). This completes the proof.
\end{proof}

This completes our characterization of graphs with two-dimensional Buchsbaum matching complexes.

\section{Concluding remarks}\label{sec:ConcludingRemarks}

Outside of dimension two, all homology manifolds that arise as matching complexes are combinatorial spheres and balls \cite{BGM}. However, in the Buchsbaum case, we do not expect higher dimensions to be as well behaved. 
Though a complete characterization is perhaps currently infeasible in general, it may be possible when restricted to certain families of graphs. For example, we can give an answer for complete bipartite graphs: If $G=K_{m,n}$ is a complete bipartite graph with $m \le n$, then $M(G)$ is Buchsbaum if and only if $n \ge 2m-2.$ This follows from \cite[Theorem~15]{Garst}, which says that $M(K_{m,n})$ is Cohen--Macaulay if and only if $n \ge 2m-1$.

We further note that \cite[Theorem~2.3]{Zi94} shows that $M(K_{m,n})$ is vertex decomposable if and only if $n \ge 2m-1$, so vertex decomposability, shellability, and Cohen--Macaulayness are equivalent for matching complexes of complete bipartite graphs. For other families of graphs, however, none of these properties hold in general. For example, if $n \ge 8$, then $M(K_n)$ is not Cohen--Macaulay (and thus not shellable or vertex decomposable). 
However, work has been done on the shellability and vertex decomposability of skeleta of these matching complexes \cite{Athanasiadis,SW07}. Thus it may be interesting to study Buchsbaumness of skeleta of matching complexes. 

Lastly, we note that some other classes of simplicial complexes can be easily described in the context of matching complexes. Recall that a complex is a matroid indepedence complex (or, simply, a matroid) if every induced subcomplex is pure. Using this definition, it is straightforward to see that $M(G)$ is a matroid if and only if $G$ does not contain a path of length three. Thus $M(G)$ is a matroid if and only if $G$ is the disjoint union of star graphs and copies of $K_3$.

\section*{Acknowledgements}\label{sec:acknowledge}

This project began as a part of the Washington Experimental Mathematics Laboratory; we thank its directors Christopher Hoffman and Jayadev Athreya. We also thank Adam Boocher and Isabella Novik for helpful discussions. B.~Goeckner was supported by a grant from the Simons Foundation (Grant Number 814268, MSRI) and by an AMS-Simons travel grant. R. Rowlands was supported by NSF grants DMS-1664865 and DMS-1953815.

\section*{Appendix}

The following is Mathematica code for computing matching complexes and checking whether a matching complex is two-dimensional and Buchsbaum. The code also checks whether a Hamiltonian graph on $7$ vertices has a Buchsbaum matching complex, which was used to generate the data in \cref{tab:c7-data}.

\begin{lstlisting}
MatchingComplex[g_] := RelationGraph[DisjointQ, EdgeList@g]
(* Computes the 1-skeleton of the matching complex. (The full matching complex is the clique complex of this graph.) *)

VertexLink[g_,v_] := VertexDelete[NeighborhoodGraph[g,v], v]
(* Computes the link of v by finding its neighborhood and deleting v *)

TwoDBuchsbaumQ[g_] := (Length@First@FindClique@g == 3) && AllTrue[VertexList@g, (ConnectedGraphQ@# && EdgeCount@#>0)& @ VertexLink[g,#]&]
(* Checks for 2D Buchsbaumness by checking that the largest clique has size 3 and all vertex links are connected graphs with at least one edge *)

c7andedges = Table[EdgeAdd[CycleGraph@7,#]& /@ Subsets[EdgeList@GraphComplement@CycleGraph@7, {i}] // DeleteDuplicates[#,IsomorphicGraphQ]&, {i,0,14}];
(* Takes all subsets of edges in the complement of C7, grouped by size of the subset, and adds them to C7, then throws out all but one graph from each isomorphism class *)

c7buchsbaums = Select[#, TwoDBuchsbaumQ@*MatchingComplex]& /@ c7andedges;
(* Picks out the graphs whose matching complex is 2D and Buchsbaum *)
\end{lstlisting}

\bibliography{bib}
\bibliographystyle{amsalpha}

\end{document}